\documentclass[11pt]{article}
\usepackage{amssymb}
\newtheorem{lem}{Lemma}[section]
\newtheorem{cor}{Corollary}[section]
\newtheorem{pro}{Proposition}[section]
\newtheorem{thm}{Theorem}[section]

\newenvironment{Literature}[1]
{\begin{thebibliography}{#1}}{\end{thebibliography}}

\begin{document}

\begin{center}
{\LARGE \bf A Local Existence Theorem  \\
\bigskip                                                                                      for the Einstein-Dirac Equation }

\bigskip  \bigskip
{\large  Eui Chul Kim}
\end{center}

\bigskip   \noindent
Department of Mathematics  \\
Inha University  \\
Inchon, 402-751, Korea    \\
e-mail: eckim@chollian.net

\bigskip \bigskip \noindent
{\bf Abstract:}   We study the Einstein-Dirac equation as well as the weak Killing equation on  Riemannian spin manifolds with codimension one foliation.
We prove that, for any manifold $M^n$ admitting real Killing spinors (resp. parallel spinors), there exist warped product metrics $\overline{\eta}$
on $M^n \times {\mathbb R}$ such that $( M^n \times {\mathbb R} , \overline{\eta} )$ admit  Einstein spinors (resp. weak Killing spinors).
To prove the result
we split the Einstein-Dirac equation into evolution equations and constraints,
by means of Cartan's frame formalism, and apply the local
preservation property of constraints.


\bigskip  \bigskip \bigskip
\noindent
\section{Introduction}

\noindent
Let $(P^m, \eta)$ be an m-dimensional smooth oriented Riemannian spin manifold and
denote by ${\rm Ric}$ and $S$ the Ricci tensor and the
scalar curvature, respectively. Let $( \, , \, ) = {\rm Re}\langle \, , \, \rangle$ be the real part of the standard Hermitian product $\langle \, , \, \rangle$ on the spinor bundle $\Sigma (P)$ over $P^m$. Let $D$ be the Dirac operator acting on sections $\psi \in \Gamma (\Sigma(P))$ of the spinor bundle $\Sigma(P)$.  The Einstein-Dirac equation is a minimal coupling of the Dirac equation
to the Einstein equation and defined by (see [11])
\[
 D \psi  =  \lambda  \psi,   \qquad
 {\rm Ric} -  \frac{1}{2}  S  \eta  =   T ,
\]
where $\lambda \in {\mathbb R}$ is some real number and the energy-momentum tensor $T$ is given by
\[
T  ( X , Y )   =   \pm \frac{1}{4}  ( X \cdot \nabla_Y \psi + Y \cdot \nabla_X  \psi   , \  \psi \, ) .
\]
A non-trivial spinor field $\psi$ solving this Einstein-Dirac system is called an {\it Einstein spinor} to eigenvalue $\lambda \in {\mathbb R}$.  In case that the scalar curvature $S$ does not vanish at any point, one defines the {\it weak Killing equation} by
\begin{eqnarray*}                                                                                                                                                                        \nabla_X \psi & =  & \frac{m}{2(m-1)S}  dS(X)  \psi  +  \frac{1}{2  (m-1)S }  X \cdot dS \cdot \psi     \cr
&     &      \cr
&     & +   \frac{2  \lambda}{ (m-2)S }  {\rm Ric}(X) \cdot \psi  -  \frac{ \lambda}{ m-2 }  X     \cdot \psi ,
\end{eqnarray*}
where $\lambda \in {\mathbb R}$ is some real number. A non-trivial solution $\psi$ to the equation is called a {\it weak Killing spinor} to weak Killing number $\lambda$ (shortly,
WK-spinor to WK-number $\lambda$). Since rescaling the length of any WK-spinor provides
an Einstein spinor, the WK-equation is stronger than the Einstein-Dirac equation (In dimension
$n=3$, the considered two equations are essentially equivalent). Moreover, the WK-equation
reduces to the Killing equation [3,8],
\[       \nabla_X \psi = - \frac{\lambda}{m} X \cdot \psi ,         \]                                                         if the metric $\eta$ is Einstein.

\bigskip
Till now, the known examples of the Einstein spinors on Riemannian manifolds are as follows: \\
\indent (i) Real Killing spinors [2,3,9,15].
\par (ii) WK-spinors on quasi-Einstein Sasakian manifolds [11].
\par (iii) Einstein spinors on product manifolds $M^6 \times N^r$, where $M^6$ is a      six-dimensional simply connected nearly K\"{a}hler manifold and $N^r$ is a
manifold of general dimension $r$ admitting Killing spinors [11].
\par (iv) WK-spinors on the three-dimensional sphere $S^3$ with non-standard merics [4,10,11].
\par (v) WK-spinors on the three-dimensional Euclidean space ${\mathbb R}^3$ with non-constant
scalar curvature [11,13].

\bigskip \noindent
The object of this paper is to establish a special existence theorem for WK- as well as Einstein spinors.
Namely, we prove the following theorem (see Theorem 5.1, Theorem 7.1 and Corollary 7.1). Interestingly, we find that
the WK-spinors constructed on ${\mathbb R}^3$ with non-constant scalar curvature (see [11], p.171) are a special case of this theorem.

\bigskip \noindent
{\bf Main Theorem:} {\it
Let $(M^n, g_M)$ be a Riemannian manifold admitting a real Killing spinor $\psi_M$. Then, for any real number $\lambda_Q \in {\mathbb R}$,
there exists a warped product metric $\overline{\eta}$ on
$Q^{n+1} = M^{n} \times {\mathbb R}$ such that $(Q^{n+1}, \overline{\eta})$ admits an Einstein spinor $\psi$ to
eigenvalue $\lambda_Q$.   In particular, if $\psi_M$ is a parallel spinor, then the
Einstein spinor $\psi$ becomes a WK-spinor to WK-number $\lambda_Q$. }

\bigskip \noindent
The key idea to prove the theorem is to split the Einstein-Dirac equation into evolution equations and constraints and apply the local
preservation property of the constraints. We will explicitly give an initial-value formulation
for the Einstein-Dirac equation, in Riemannian setting,
  and solve it for a specific class of initial
data sets.  It is well-known that, in Riemannian signature
the Einstein equations are generally of elliptic type, making the
initial-value problem (the Cauchy problem) for general smooth data
inappropriate. However, when the considered Riemannian manofolds
admit a codimension one foliation, one can represent the Einstein
equations to be of hyperbolic type, just as one does over
Lorentzian manifolds, and can indeed formulate the initial-value problem in a natural way.
\par
So far, not much has been studied about the initial-value problem for the
Einstein-Dirac equation. In Lorentzian signature, the spacelike
initial-value problem for the Einstein-Dirac system was considered
by Bao/Isenberg/Yasskin [1] in terms of 3+1 Hamiltonian formalism,
but no existence theorem was proved there. Recently,
Friedrich/Rendall indicated [7], in terms of Penrose's two-spinor
formalism, that the Einstein-Dirac equation may be reduced to
symmetric hyperbolic evolution equations, illustrating some
questions arising in the reduction.
\par
In this paper we give an invariant description of the initial-value formulation for the       Einstein-Dirac equation on Riemannian manifolds with codimension one foliation, in an explicit form and in complete generality.  The splitting of the Einstein-Dirac equation into
evolution equations and constraints will be achieved in terms of Cartan's frame formalism, and hence our formulation is valid
on Riemannian manifolds $M^n \times {\mathbb R}$ of general dimension $n+1$.
The first three sections (Section 2,3,4) of the paper are devoted to establishing the
basic framework, the {\it hyperbolic representation} of curvatures and the Dirac equation, on (possibly compact) manifolds with codimension one foliation, and the framework
may be of independent interest for further study of the behaviour of spinor field equations
under global change of metrics.

\bigskip  \noindent
\section{Representation of curvatures and the Dirac equation with respect to reference metric}

\noindent
Let $P^m$ be an m-dimensional simply-connected smooth oriented manifold allowing spin structure, and let $\eta, \overline{\eta}$ be
two Riemannian metrics on $P^m$. Henceforth we fix the notation $\eta$ to denote a reference metric.   Then there exists a
unique $(1,1)$-tensor field $K$ on $P^m$ that is positive definite with respect to $\eta$ and
satisfies
\[         \overline{\eta} (X,Y ) = \eta ( K(X), K(Y) )           \]
for all vector fields $X, Y$.
Recall that
the Levi-Civita connection $\nabla^{\overline{\eta}}$ of $(P^m, \overline{\eta} )$ is characterized by the Koszul formula
\begin{eqnarray*}
2 \, \overline{\eta} ( Z , \, \nabla^{\overline{\eta}}_X Y ) & = & X \{ \overline{\eta} ( Y , Z) \}
+  Y \{ \overline{\eta} ( Z , X) \}  -  Z \{ \overline{\eta} (X , Y) \}    \\
&      &     \\                  &    &   +
 \overline{\eta} ( Z  ,  [ X , Y ] )  +   \overline{\eta} ( Y  ,  [ Z , X ] )  - \overline{\eta} ( X,  [ Y , Z ] )  .
\end{eqnarray*}
Letting $(E_1,   \ldots, E_m )$ be a local $\eta$-orthonormal frame field on $P^m$, for which  \[         \Big( F_1 = K^{-1}(E_1),  \ldots, F_m = K^{-1}(E_m) \Big)               \]                                 is $\overline{\eta}$-orthonormal, and
inserting $X = F_i, \, Y = F_j,  \, Z = F_k$ into the Koszul formula, we have
\begin{eqnarray*}
&     & 2 \, \eta \Big(  E_k , \, K \{ \nabla^{\overline{\eta}}_{K^{-1}(E_i)} (
K^{-1} E_j)  \}   \Big)     \\   &     &     \\  & = &  \eta \Big( E_k , \,
K \{ [  K^{-1}(E_i) ,  K^{-1}(E_j) ] \} \Big)
 +  \eta \Big( E_j ,  K \{ [  K^{-1}(E_k) ,  K^{-1}(E_i) ] \} \Big)         \\
&     &     \\    &      & -  \eta \Big( E_i ,  K \{ [  K^{-1}(E_j) ,
K^{-1}(E_k) ] \} \Big)       \\   &     &     \cr &  =  &   \eta \Big(  E_k
 , \, K(\nabla^{\eta}_{F_i}   F_j)  -
K(\nabla^{\eta}_{F_j} F_i)  \Big)
 +   \eta \Big(  E_j  , \, K(\nabla^{\eta}_{F_k} F_i)  -  K(\nabla^{\eta}_{F_i} F_k)  \Big)          \\
&     &     \\
&     & -    \eta \Big( E_i  , \, K(\nabla^{\eta}_{F_j}
F_k)  -  K(\nabla^{\eta}_{F_k}
F_j)  \Big)    \\  &    &       \\  &  =  &   \eta \Big( E_k
, \, K\{ (\nabla^{\eta}_{F_i} K^{-1})(E_j) \} +
\nabla^{\eta}_{F_i} E_j  -   K \{
(\nabla^{\eta}_{F_j} K^{-1}) (E_i) \}  -
\nabla^{\eta}_{F_j} E_i   \Big)    \\   &    &      \\   &
& +   \eta \Big( E_j , \, K \{ (\nabla^{\eta}_{F_k}
K^{-1})(E_i) \}  +  \nabla^{\eta}_{F_k} E_i  -   K \{
(\nabla^{\eta}_{F_i} K^{-1}) (E_k) \}  -
\nabla^{\eta}_{F_i} E_k  \Big)    \\   &    &      \\   &    &
-   \eta \Big( E_i  , \, K \{ (\nabla^{\eta}_{F_j}
K^{-1})(E_k) \} +  \nabla^{\eta}_{F_j} E_k  -   K \{
(\nabla^{\eta}_{F_k} K^{-1}) (E_j) \} -
\nabla^{\eta}_{F_k} E_j  \Big)    \\   &     &     \\   & = &
2 \,  \eta \Big( E_k   , \,  \nabla^{\eta}_{K^{-1}(E_i)} \, E_j  \Big)
\cr &     &     \\   &     & +  \eta  \Big( E_k  , \, K \{  (
\nabla^{\eta}_{K^{-1}(E_i)} K^{-1} ) (E_j) \}  - K \{  (
\nabla^{\eta}_{K^{-1}(E_j)}  K^{-1} )(E_i) \}  \Big)    \\   &     &
\\   &     & +  \eta  \Big( E_j  , \, K \{  (
\nabla^{\eta}_{K^{-1}(E_k)}  K^{-1} )(E_i) \}  - K \{  (
\nabla^{\eta}_{K^{-1}(E_i)}  K^{-1} )(E_k) \}  \Big)    \\    &     &
\\  &    & -   \eta  \Big( E_i  , \, K \{  (
\nabla^{\eta}_{K^{-1}(E_j)}  K^{-1} )(E_k) \}  - K \{  (
\nabla^{\eta}_{K^{-1}(E_k)}  K^{-1} )(E_j) \}  \Big)    .
\end{eqnarray*}
Thus we obtain the following formula.

\bigskip \noindent
\begin{pro}                                                                                   The Levi-Civita connections $\nabla^{\overline{\eta}}, \,
\nabla^{\eta}$ are related by
\[
\nabla^{\overline{\eta}}_{K^{-1}(X)} \left( K^{-1}(Y) \right)  =
K^{-1} \left( \nabla^{\eta}_{K^{-1}(X)} Y \right) +  K^{-1} \left\{ \Lambda_{\eta} (X, Y) \right\} ,
\]
where $ \Lambda_{\eta}$ is the (1,2)-tensor field defined by
\begin{eqnarray*}
2 \, \eta ( \Lambda_{\eta} (X , Y ),  Z)
&  =  &   \eta \left(  Z  , \, K \{  ( \nabla^{\eta}_{K^{-1}(X)} K^{-1} ) (Y) \}  -
K \{  ( \nabla^{\eta}_{K^{-1}(Y)}  K^{-1} )(X) \}  \right)    \\
&     &     \\
&    & +  \eta \left(  Y  , \, K \{  ( \nabla^{\eta}_{K^{-1}(Z)}  K^{-1} )(X) \}  -
K \{  ( \nabla^{\eta}_{K^{-1}(X)}  K^{-1} )(Z) \}  \right)    \cr
&     &     \\
&     & +  \eta \left(  X  , \, K \{  ( \nabla^{\eta}_{K^{-1}(Z)}  K^{-1} )(Y) \}  -
K \{  ( \nabla^{\eta}_{K^{-1}(Y)}  K^{-1} )(Z) \} \right)   .
\end{eqnarray*}
\end{pro}

\bigskip \noindent
{\bf Remark 2.1} (i) The exact difference between the Levi-Civita connections $\nabla^{\overline{\eta}},
\, \nabla^{\eta}$ is related to the tensor $\Lambda_{\eta}$ by
\[         \nabla^{\overline{\eta}}_X Y - \nabla^{\eta}_X Y  =
K^{-1} \{ \Lambda_{\eta} (KX, KY) \} + K^{-1} \{ ( \nabla^{\eta}_X K)(Y) \} .         \]
\indent (ii)  The relation
\[             \eta( \Lambda_{\eta} ( X, Z ) , Y )
+ \eta( \Lambda_{\eta} ( X , Y ) , Z ) = 0
\]                is valid for all vector fields $X, Y, Z$.   \\
\indent (iii)  Since    \[         \Lambda_{\eta} (X, Y) -
\Lambda_{\eta} (Y,X) = K \left\{ ( \nabla^{\eta}_{K^{-1}X} K^{-1} )
(Y)  \right\} -
 K \left\{ ( \nabla^{\eta}_{K^{-1}Y} K^{-1} ) (X)  \right\}  ,
\]
$\Lambda_{\eta} \equiv 0$ vanishes identically if and only if $
\Lambda_{\eta} (X, Y) = \Lambda_{\eta} ( Y , X)$ for all vector fields $X, Y$.

\bigskip \noindent
We will often use the shorthand notation $\Lambda_{\eta} = \Lambda$, if there is no
possibility of confusion. Proposition 2.1 enables us to describe the behaviour of
curvatures under global change of metrics in a nice way : A direct computation gives
\begin{eqnarray*}
&     & R_{\overline{\eta}} ( K^{-1}X , K^{-1} Z ) ( K^{-1} Y)  -
K^{-1} \left\{ R_{\eta}
(K^{-1}X , K^{-1} Z ) (Y) \right\}     \\
&     &    \\
&   = &   K^{-1} \left\{ ( \nabla^{\eta}_{K^{-1}(X)} \Lambda ) (Z, Y) -
( \nabla^{\eta}_{K^{-1}(Z)} \Lambda ) (X, Y) \right\}     \\ &     & \\
&     &  +  K^{-1} \Big\{ \Lambda ( X , \Lambda (Z, Y ) ) -
\Lambda ( Z , \Lambda (X, Y ) ) \Big\}     \\
&    &      \\
&    &  + K^{-1} \Big\{ \Lambda ( \Lambda (Z, X) - \Lambda (X, Z), \, Y ) \Big\},
\end{eqnarray*}
where  $R_{\overline{\eta}}$ (resp. $R_{\eta}$)
is the Riemann tensor of $\overline{\eta}$ (resp. $\eta$).
Contracting both sides of the equation, we can now represent the Ricci curvature ${\rm Ric}_{\overline{\eta}}$ as well as the scalar curvature $S_{\overline{\eta}}$ with respect to the reference metric $\eta$.

\bigskip \noindent
\begin{pro}
\begin{eqnarray*}
&     & {\rm Ric}_{\overline{\eta}} (K^{-1}Y, K^{-1}Z) -
\sum_{j=1}^m \eta \Big( E_j, \ R_{\eta} ( K^{-1} E_j, K^{-1}Z)(Y) \Big)     \\
&     &     \\
& = & \sum_{j=1}^m  \eta \Big( E_j, \ ( \nabla^{\eta}_{K^{-1}E_j} \Lambda ) (Z, Y)
- ( \nabla^{\eta}_{K^{-1}Z} \Lambda ) (E_j, Y)  \Big)      \\
&      &     \\
&     & + \sum_{j=1}^m  \eta \Big( E_j, \ \Lambda (E_j, \Lambda(Z, Y))
- \Lambda ( \Lambda(E_j, Z), Y)   \Big) ,
\end{eqnarray*}
In particular,
\begin{eqnarray*}
&     & S_{\overline{\eta}} - \sum_{i,j =1}^m \eta \Big( E_i , R_{\eta}( K^{-1}E_i ,
K^{-1} E_j ) (E_j) \Big)   \\ &     &     \\    & = & 2 \sum_{i,j=1}^m
\eta \Big( E_i , ( \nabla^{\eta}_{K^{-1} (E_i)} \Lambda ) ( E_j , E_j )
\Big)      - \sum_{i,j,k=1}^m \eta( E_k , \Lambda ( E_i , E_i )
) \eta( E_k , \Lambda( E_j , E_j ) )       \\    &         &    \\
&     & -  \sum_{i,j,k=1}^m \eta( E_k , \Lambda ( E_i , E_j
) ) \eta( E_k , \Lambda ( E_j , E_i ) )  .
\end{eqnarray*}
\end{pro}

\bigskip
Next, we review briefly the behaviour of the Dirac operator under
change of metrics. Let $T(P)$ be the tangent bundle of $P^m$, and
let $\Sigma(P)_{\overline{\eta}}$ (resp. $\Sigma(P)_{\eta}$) be the spinor bundle of
$(P, \overline{\eta})$ (resp. $(P, \eta)$) equipped with the standard Hermitian product $\langle \, , \, \rangle_{\overline{\eta}}$ (resp.
$\langle \, , \, \rangle_{\eta}$).
We know that there exists a  natural isomorphism                                                                                                              $ {\widetilde{K}} : \Sigma {(P)}_{\overline{\eta}} \longrightarrow \Sigma {(P)}_{\eta} $ with
\[
 \langle {\widetilde{K}} ( \varphi) \, , \, {\widetilde{K}} ( \psi) \rangle_{\eta} \, = \, \langle \varphi , \, \psi  \rangle_{\overline{\eta}} ,   \qquad
 ( K Z ) \cdot ( {\widetilde{K}} \psi ) \, = \, {\widetilde{K}} ( Z \cdot \psi )   ,
\]
for all $Z \in T(P)  , \  \varphi , \psi \in \Sigma {(P)}_{\overline{\eta}} $, where
 the dots "$\cdot$" in the latter relation indicate the Clifford multiplication with respect to $\eta$ and $\overline{\eta}$, respectively.
In terms of local $\overline{\eta}$-orthonormal frame field $(F_1, \ldots, F_m)$, the spin derivative $\nabla^{\overline{\eta}} \varphi$ is  expressed as
\[             \nabla^{\overline{\eta}}_X \varphi = X(\varphi) + \frac{1}{4} \sum_{i=1}^m
F_i \cdot \nabla^{\overline{\eta}}_X F_i \cdot \varphi ,  \qquad  \varphi  \in \Gamma (\Sigma (P)_{\overline{\eta}} ) ,       \]
and the Dirac operator $D_{\overline{\eta}} \varphi$ as
\[     D_{\overline{\eta}} \varphi = \sum_{j=1}^m F_j \cdot \nabla^{\overline{\eta}}_{F_j} \varphi.
\]
Making use of the formula in Proposition 2.1,
one finds now readily that the  spinor derivatives $\nabla^{\overline{\eta}} , \, \nabla^{\eta}$ and the Dirac operators $D_{\overline{\eta}}, \, D_{\eta}$ are related as follows.

\bigskip \noindent
\begin{pro}                                                                                   (see [5]) \ For all \ $\psi \in \Gamma (\Sigma (P)_{\eta})$,
\begin{eqnarray*}
\Big\{  \widetilde{K} \circ  \nabla^{\overline{\eta}}_{K^{-1}(E_j)}  \circ  ( \widetilde{K} )^{-1}   \Big\} (\psi)
&  =  &  \nabla^{\eta}_{K^{-1}(E_j)} \psi + \frac{1}{4} \sum_{k,l=1}^m  \Lambda_{j k l} E_k \cdot E_l \cdot \psi    ,    \\
&     &    \\
&     &    \\
 \Big\{  \widetilde{K} \circ D_{\overline{\eta}} \circ (\widetilde{K} )^{-1}  \Big\} ( \psi )
&  =  &  \sum_{i=1}^m   E_i \cdot \nabla^{\eta}_{K^{-1}(E_i)} \psi + \frac{1}{4}
\sum_{j,k,l=1}^m   \Lambda_{jkl} E_j \cdot E_k \cdot E_l \cdot \psi        \\
&      &       \\
& = &  \sum_{i=1}^m  E_i \cdot \nabla^{\eta}_{K^{-1}(E_i)} \psi
- \frac{1}{2} \sum_{j,k=1}^m \Lambda_{jjk} E_k \cdot \psi           \\
&     &       \\
&      & \quad + \frac{1}{2}  \sum_{j<k<l}^m   ( \Lambda_{jkl} + \Lambda_{klj} + \Lambda_{ljk} ) E_j \cdot E_k \cdot E_l \cdot \psi ,
\end{eqnarray*}
where $\Lambda_{jkl} : = \eta ( \Lambda_{\eta}(E_j, E_k) , E_l )$.
\end{pro}

\bigskip \noindent
\section{Representation of curvatures on manifolds with codimension one foliation}

\noindent
In this section we establish an intrinsic setting of the formulas that constitute the well-known evolution system
for the Einstein (vacuum) equation (see [6,7]).  The evolution system consists of
 two differential equations, describing the evolution of metrics (see Corollary 3.1) and the evolution of symmetric
(0,2)-tensor fields (see Proposition 3.2), respectively.
The main aim of this section is to represent the Ricci tensor ${\rm Ric}_{\overline{\eta}}$
{\it hyperbolically} with respect to codimension one foliation.  We
use the terminology "{\it hyperbolic} representation" in the sense that such representation
of differential operators, on manifolds with codimension one foliation, transforms field equations of elliptic type involving metrics to hyperbolic systems
in PDE theory. Note in this view that the formulas in Proposition 2.2 may be
thought of as the {\it elliptic} representation of curvatures.
\par
Let $(Q^{n+1}, \eta)$ be an (n+1)-dimensional smooth oriented Riemannian spin manifold. We assume that there exists a codimension one foliation on $(Q^{n+1}, \eta)$ defined by a unit vector field $E_{n+1}$
with $dE^{n+1} = 0$, where $E^{n+1} = \eta ( E_{n+1}, \cdot )$ is the dual
1-form of $E_{n+1}$.                                                                          Letting $E_{n+1}^{\bot}$ denote the $\eta$-orthogonal
complement  of $E_{n+1}$ in the tangent bundle $T(Q)$,
we note that $dE^{n+1} = 0$ implies the following facts(e.g. see [14]):
\par (i)   For all vector fields $V, W$ belonging to $E_{n+1}^{\bot}$,  all of $[V, W], \, \nabla^{\eta}_{E_{n+1}} V$ and $\nabla^{\eta}_W E_{n+1}$ belong to
 $E_{n+1}^{\bot}$  .
\par (ii) $\nabla^{\eta}_{E_{n+1}} E_{n+1} = 0$ .
\par (iii) If $Q^{n+1}$ is compact, then all the slices of the foliation are diffeomorphic.
\par (iv) If $Q^{n+1}$ is simply-connected, then $E^{n+1} = ds$ for some real-valued function $s : Q^{n+1} \longrightarrow {\mathbb R}$ ($Q^{n+1}$ must be noncompact) and the foliation is defined by the level hypersurfaces $s = \mbox{\it constant}$.

\bigskip
Let $(E_1, \ldots, E_n, E_{n+1})$
be a local $\eta$-orthonormal frame field on $Q^{n+1}$, with $E_j \in E_{n+1}^{\bot}, j=1, \cdots, n$,
and  $(E^1, \ldots, E^n, E^{n+1})$ the dual frame field.
                                      Denote
$\otimes^r_s(E_{n+1}^{\bot})$ the space of all $(r,s)$-tensor fields
$B$ on $Q^{n+1}$ such that        \[              \eta ( E_{i_1} \otimes \cdots \otimes E_{i_r}, \
B( E_{j_1} \otimes \cdots \otimes E_{j_s}) ) = 0 ,      \]                 whenever either $i_k = n+1$
for some $i_k$ or $j_l = n+1$ for some $j_l$.
Now, consider a positive definite (1,1)-tensor field $K$ on $(Q^{n+1}, \eta)$. Letting
$\overline{\eta}$ be the metric induced by $K$ via $\overline{\eta} (X, Y) = \eta (K(X), K(Y) )$ and identifying $\overline{\eta}$ with $K^2$,
we can express $\overline{\eta}$ as
\begin{eqnarray*}
\overline{\eta} = K^2 & = & \Big\{ \sum_{i,j =1}^n (L^2)_j^i E^j \otimes E_i  \Big\}
+  E^{n+1} \otimes L^2(\zeta) +  \eta(
L^2(\zeta), \, \cdot \, ) \otimes E_{n+1}       \cr &     &
\cr &     &  \quad + \Big\{ \eta( L(\zeta), L(\zeta) ) +
 \rho^2 \Big\} E^{n+1} \otimes E_{n+1} ,
\end{eqnarray*}
where $L \in  \otimes^1_1(E_{n+1}^{\bot}), \ \zeta \in
\otimes^1_0(E_{n+1}^{\bot})$ and $\rho : Q^{n+1} \longrightarrow  {\mathbb
R}$ is a positive function.  \\
This may be thought of as an intrinsic (Riemannian) version of the well-known ADM-representation of metrics in general relativity.
$\zeta$ agrees with the {\it shift vector field} and $\rho$ with the {\it lapse function}.
Note that the (1,1)-tensor $K^2$ is related to the (1,1)-tensor $L^2$ by
\[
K^2(V) =  L^2(V) + \eta( L^2(\zeta), V) E_{n+1}
\]
for all vector fields $V \in \otimes^1_0(E_{n+1}^{\bot})$ and
\[            K^2(E_{n+1})  =  L^2(\zeta) +   \{
 \eta( L(\zeta), L(\zeta) ) +  \rho^2 \}  E_{n+1}.         \]
Furthermore,
$L^2$ is positive definite on each slice of the foliation and, on the slices, coincides with the
metrics induced by $K^2$.
Certainly,
\[
\Big( F_1 :=L^{-1}(E_1), \ldots, F_n :=L^{-1}(E_n),  F_{n+1} :=
\rho^{-1} (E_{n+1} - \zeta ) \Big)       \]                   is a local $\overline{\eta}$-orthonormal frame field on $Q^{n+1}$, its dual frame field being
given by       \[          F^i = L(E^i) + \eta( L(\zeta), E_i )
E^{n+1},  \qquad F^{n+1} = \rho E^{n+1} .         \]

\par
Let $Z$ be a vector field on $Q^{n+1}$, and let
$V, W$ be vector fields in $E_{n+1}^{\bot}$.  In what follows we fix the notations $V, W$
to mean vector fields in $E_{n+1}^{\bot}$.
 Then, one verifies easily the following basic identities:
\begin{eqnarray*}
\overline{\eta}(V, W) & = & \eta( L^2(V), W) = \eta (V, L^2(W) ) =
\eta (L(V), L(W) )  ,      \\   &     &      \\     \overline{\eta}(V,
E_{n+1}) & = & \eta( V, L^2(\zeta) ) ,       \\   &      &      \\
\overline{\eta}(E_{n+1}, E_{n+1}) & = & \eta( L(\zeta), L(\zeta) )
+  \rho^2,       \\    &     &     \\
\overline{\eta}(Z, F_{n+1}) & = &  \rho  \eta( Z, E_{n+1}) ,
\\    &     &     \\     \overline{\eta}(V, F_{n+1}) & = & 0.
\end{eqnarray*}
The identity $\overline{\eta}(V, F_{n+1})  =  0$ in the last line
implies that $E_{n+1}^{\bot}$ coincides with the $\overline{\eta}$-orthogonal
complement of $F_{n+1}$ in $T(Q)$.

\par
We let           \[ {\rm II} (V) : = - \nabla^{\overline{\eta}}_V F_{n+1}    \qquad   \mbox{and}
\qquad
\Theta(V) : = - \nabla^{\eta}_V E_{n+1}       \]           denote the second fundamental form,
on each slice, defined by the unit vector field $F_{n+1}$ and $E_{n+1}$,
respectively. Let
$\overline{g}$ (resp. $g$) denote the metric, on each slice,
 induced by $\overline{\eta}$ (resp. $\eta$) and
 $\nabla^{\overline{g}}$ (resp. $\nabla^g$) its Levi-Civita
connection. In the notations, the tensor $L$ satisfies
\begin{eqnarray*}
\nabla^{\eta}_{E_{n+1}} L  & \in  & \otimes^1_1 ( E_{n+1}^{\bot} )   ,   \\
&     &     \\
( \nabla^{\eta}_V L ) (W) & = & ( \nabla^g_V L )(W) + \Theta(V, L(W)) E_{n+1}  ,   \\
&     &     \\
( \nabla^{\eta}_V L)(E_{n+1}) & = & L( \Theta(V) ) .
\end{eqnarray*}

\par
In order to represent curvatures of $\overline{\eta}$ hyperbolically  with respect to codimension one foliation (Proposition 3.2), we must explicitly
know how the connection $\nabla^{\overline{\eta}}$ is related to
the connections $\nabla^{\eta}$ and $\nabla^g$. This is done in the following proposition, which may be thought of as the hyperbolic version of Proposition 2.1.

\bigskip  \noindent
\begin{pro}
\begin{eqnarray*}
&  (i)   & \quad \nabla^{\overline{\eta}}_{F_i} F_j = \nabla^{\overline{g}}_{F_i} F_j
+ {\rm II} (F_i, F_j) F_{n+1}    \quad (1 \leq i, j \leq n)     \\
&      &     \\
&     & = L^{-1} (\nabla^g_{L^{-1}E_i} E_j) + L^{-1} \{ \Lambda_g (E_i, E_j) \} +
{\rm II} (L^{-1}E_i, L^{-1} E_j) \{ \rho^{-1} (E_{n+1} - \zeta ) \},     \\
&     &     \\
&      &     \\
&  (ii)    &  \quad  \rho \cdot  \nabla^{\overline{\eta}}_{F_{n+1}} F_j
        \\    &     &      \\    &    & = L^{-1}(
\nabla^{\eta}_{E_{n+1}} E_j ) - L^{-1} (\nabla^g_{\zeta} E_j ) +
\frac{1}{2} ( \nabla^{\eta}_{E_{n+1}} L^{-1}) (E_j) - \frac{1}{2}
( \nabla^g_{\zeta} L^{-1} ) (E_j)        \\     &      &     \\     &
& \quad + \frac{1}{2} \nabla^g_{L^{-1}E_j} \zeta + \frac{1}{2} \Theta
(L^{-1} E_j) + \frac{1}{2} \sum_{i=1}^n \eta \Big( E_j, \ (
\nabla^{\eta}_{E_{n+1}} L ) (L^{-1} E_i)  \Big) L^{-1} E_i        \\
&      &     \\    &      & \quad - \frac{1}{2} \sum_{i=1}^n \eta\Big( E_j, \
(\nabla^g_{\zeta} L)(L^{-1}E_i) + L( \nabla^g_{L^{-1} E_i} \zeta )
+ (L \circ \Theta \circ L^{-1})(E_i)  \Big) L^{-1}E_i         \\    &
&       \\
&      & \quad + d \rho (L^{-1} E_j) F_{n+1} ,    \\    &      &      \\
&      &     \\    &  (iii)    &  \quad   - \rho \cdot
\nabla^{\overline{\eta}}_{F_j} F_{n+1} =  \rho \cdot {\rm II}
(L^{-1} E_j)         \\    &     &      \\     &     &  =  \frac{1}{2} (
\nabla^{\eta}_{E_{n+1}} L^{-1} )(E_j) - \frac{1}{2} \sum_{i=1}^n
\eta\Big( E_j, \ ( \nabla^{\eta}_{E_{n+1}} L ) (L^{-1} E_i)  \Big) L^{-1}
E_i        \\    &      &     \\    &      &  \quad - \frac{1}{2} (
\nabla^g_{\zeta} L^{-1} ) (E_j) + \frac{1}{2} \nabla^g_{L^{-1}E_j}
\zeta + \frac{1}{2} \Theta (L^{-1} E_j)      \\    &     & \quad +
\frac{1}{2} \sum_{i=1}^n \eta\Big( E_j, \ (\nabla^g_{\zeta}
L)(L^{-1}E_i) + L( \nabla^g_{L^{-1} E_i} \zeta ) + (L \circ \Theta
\circ L^{-1})(E_i)  \Big) L^{-1}E_i  ,       \\    &      &      \\    &
&     \\     &  (iv)    &  \nabla^{\overline{\eta}}_{F_{n+1}} F_{n+1} = -
\rho^{-1} \sum_{i=1}^n d \rho (L^{-1} E_i) L^{-1} E_i.
\end{eqnarray*}
\end{pro}

\noindent
{\bf Proof.} One computes directly, substituting the identities
\begin{eqnarray*}
&      &  [
F_{n+1}, V ]  =   \nabla^{\overline{\eta}}_{F_{n+1}} V + {\rm II}(V)    \\
&     &    \\
&  =  & \rho^{-1} d \rho (V) F_{n+1} + \rho^{-1} [ E_{n+1}, V ] +
\rho^{-1} [ V, \zeta ]    \\
&     &    \\
&  =  &  \rho^{-1} d\rho(V) F_{n+1} + \rho^{-1}
\nabla^{\eta}_{E_{n+1}} V + \rho^{-1} \Theta(V) + \rho^{-1} [V,
\zeta ].
\end{eqnarray*}
in the Koszul formula. \quad $\rule{2mm}{2mm}$

\bigskip  \noindent
We can equivalently rewrite the third equation (iii) in Proposition 3.1 as follows.

\bigskip  \noindent
\begin{cor}
\begin{eqnarray*}
2 \rho \cdot {\rm II} (V, W)
& = & - \eta ( ( \nabla^{\eta}_{E_{n+1}} L)(V), \ L(W) ) -
\eta ( L(V), \  ( \nabla^{\eta}_{E_{n+1}} L)(W) )     \\     &     &
\\   &    & + \eta ( ( \nabla^g_{\zeta} L)(V), \ L(W) ) + \eta ( L(V),
\  ( \nabla^g_{\zeta} L)(W) )     \\   &     &     \\   &     & +
\eta ( \nabla^g_V \zeta + \Theta(V) , \ L^2(W) )
+
 \eta ( L^2(V), \  \nabla^g_W \zeta + \Theta(W) ) .
\end{eqnarray*}
\end{cor}

\bigskip
To prove Proposition 3.2 below, we need to recall the Gauss-Codazzi equations that relates the curvatures of $(Q^{n+1}, \overline{\eta} )$ to the curvatures of the slices.

\bigskip  \noindent
\begin{lem}
\begin{eqnarray*}
&  (i)   &  \quad  R_{\overline{\eta}}(V_1,V_2)(V_3) - R_{\overline{g}}  (V_1,V_2)(V_3)
\\   &      &     \\   &     &  =  {\rm II} (V_1, V_3){\rm
II}(V_2) - {\rm II}(V_2, V_3) {\rm II}(V_1)    +  \Big\{ ( \nabla^{\overline{g}}_{V_1} {\rm II}
)(V_2, V_3) - ( \nabla^{\overline{g}}_{V_2} {\rm II} )(V_1, V_3) \Big\} F_{n+1},
\\    &     &     \\   &     &     \\   &   (ii)    &  \quad {\rm
Ric}_{\overline{\eta}}(W) -  R_{\overline{\eta}}(W,
F_{n+1}) (F_{n+1}) - {\rm Ric}_{\overline{g}}(W)        \\   &
&     \\    &     &  =  ( {\rm II} \circ {\rm II})(W) -
 ( {\rm Tr}_{\overline{g}} {\rm II}) {\rm II}(W)  +  \Big\{ d ( {\rm Tr}_{\overline{g}}
{\rm II}) (W) - ({\rm div}_{\overline{g}}  {\rm II}) (W)  \Big\}
F_{n+1},   \\   &     &     \\   &     &     \\   &  (iii)    &   \quad
S_{\overline{\eta}} - 2 \cdot {\rm
Ric}_{\overline{\eta}}(F_{n+1}, F_{n+1})  - S_{\overline{g}}  =   {\rm
Tr}_{\overline{g}}( {\rm II}^2 ) -  ( {\rm Tr}_{\overline{g}} {\rm
II})^2 .
\end{eqnarray*}
\end{lem}

\bigskip  \noindent
\begin{pro}
\begin{eqnarray*}
 &       & ( \nabla^{\eta}_{E_{n+1}} {\rm II}) (V, W)     \\
&       &      \\    &  =  &  \rho  \cdot  \Big\{   {\rm
Ric}_{\overline{\eta}}(V,W) -  {\rm
Ric}_{\overline{g}}(V,W) - 2 \overline{\eta} (V,  {\rm II}^2(W) )
+  ( {\rm Tr}_{\overline{g}} {\rm II} ) {\rm II}(V, W)   \Big\}
\\   &       &     \\   &       & + ( \nabla^{\overline{g}}_{\zeta}
{\rm II} ) (V,W) +  \overline{\eta} ( {\rm II}(V),
\nabla^{\overline{g}}_W \zeta ) +   \overline{\eta} ( {\rm II}(W),
\nabla^{\overline{g}}_V \zeta )         \\    &       &     \\    &
& + \overline{\eta} ( \Theta(V),  {\rm II}(W) )  +
\overline{\eta} ( \Theta(W),  {\rm II}(V) ) +  (
\nabla^{\overline{g}}_W d \rho ) (V) .
\end{eqnarray*}
\end{pro}

\noindent
{\bf Proof.} Via a direct computation, we have
\begin{eqnarray*}
&     & R_{\overline{\eta}} (W, F_{n+1})(F_{n+1})  =
\nabla^{\overline{g}}_W (\nabla^{\overline{\eta}}_{F_{n+1}}
F_{n+1}) + \rho^{-1} d \rho(W) \nabla^{\overline{\eta}}_{F_{n+1}}
F_{n+1}    - ({\rm II} \circ {\rm II})(W)        \\   &      &
\\   &      & \qquad  \quad  + \rho^{-1} [ E_{n+1},  {\rm II}(W) ]
- \rho^{-1} [ \zeta,  {\rm II}(W) ] - \rho^{-1} {\rm II}( [
E_{n+1},  W] ) - \rho^{-1} {\rm II} ([W, \zeta ]).
\end{eqnarray*}
Using the equation (iv) in Proposition 3.1, we compute
\begin{eqnarray*}
&       &  \overline{\eta} (V, \  R_{\overline{\eta}} ( W, F_{n+1}
) ( F_{n+1}) )        \\   &      &      \\   &  =  &  W \{
\overline{\eta} ( V,  \nabla^{\overline{\eta}}_{F_{n+1}} F_{n+1} )
\} -  \overline{\eta} ( \nabla^{\overline{g}}_W V,
\nabla^{\overline{\eta}}_{F_{n+1}} F_{n+1} ) + \rho^{-1} d \rho(W)
\overline{\eta} ( V,  \nabla^{\overline{\eta}}_{F_{n+1}} F_{n+1} )
\\   &      &      \\    &     & - \overline{\eta} (V,   {\rm II}^2
(W) ) + \rho^{-1}  \overline{\eta} (V,  [ E_{n+1},  {\rm II}(W) ]
) -  \rho^{-1}  \overline{\eta} (V,  [ \zeta,  {\rm II}(W) ] )
\\    &      &       \\   &     & - \rho^{-1} \overline{\eta} ( {\rm
II}(V) , [ E_{n+1}, W] ) +
 \rho^{-1} \overline{\eta} ( {\rm II}(V) , [ \zeta, W] )         \\
&     &      \\   &  =  &   \rho^{-2} d \rho(V) d \rho(W)
-  \rho^{-1} ( \nabla^{\overline{g}}_W d \rho )(V)
 -  \rho^{-1} d \rho (  \nabla^{\overline{g}}_W  V )          \\
&     &      \\   &     & +  \rho^{-1} d \rho (
\nabla^{\overline{g}}_W  V )    -   \rho^{-2} d \rho(V) d
\rho(W) - \overline{\eta} ( V,  {\rm II}^2 (W) )         \\   &
&       \\   &     & +  \overline{\eta} \Big( V, \, [ F_{n+1},  {\rm
II}(W) ] - \rho^{-1} d \rho( {\rm II}(W) ) F_{n+1}  \Big) -
\overline{\eta} \Big( {\rm II}(V),  \, [ F_{n+1},  W ] - \rho^{-1} d
\rho( W ) F_{n+1}  \Big)          \\   &      &      \\   &  =  &  -
 \rho^{-1} ( \nabla^{\overline{g}}_W d \rho )(V) -
\overline{\eta} ( V,  {\rm II}^2 (W) )   +  \overline{\eta} ( V,  [ F_{n+1},  {\rm II}(W) ]  ) -
\overline{\eta} ( {\rm II}(V),  [ F_{n+1},  W ]  ).
\end{eqnarray*}

\noindent On the other hand, from the Koszul formula for $\nabla^{\overline{\eta}}$, we know that
\begin{eqnarray*}
&     & 2  \overline{\eta} (V,  \nabla^{\overline{\eta}}_{F_{n+1}}
\{ {\rm II}(W)  \}  )       \\   &     &     \\   & =  & F_{n+1}  \{
\overline{\eta} (V,  {\rm II}(W) )  \}  +  \overline{\eta} ( V,  [
F_{n+1},  {\rm II}(W) ]  ) -  \overline{\eta} ( {\rm II}(W),  [
F_{n+1},  V ]  )         \\   &     &     \\   &  =  & 2  F_{n+1}
\{  \overline{\eta} (V,  {\rm II}(W) )  \}  + 2  \overline{\eta} (
\nabla^{\overline{\eta}}_{F_{n+1}} V, \nabla^{\overline{\eta}}_W
F_{n+1}  ),
\end{eqnarray*}
which gives
\begin{eqnarray*}
&      &  \overline{\eta} ( V,  [ F_{n+1},  {\rm II}(W) ]  ) -
\overline{\eta} ( {\rm II}(W),  [ F_{n+1},  V ]  )         \\   &
&     \\   &  =  &  F_{n+1}  \{  \overline{\eta} (V,  {\rm II}(W) )
\}  + 2  \overline{\eta} ( \nabla^{\overline{\eta}}_{F_{n+1}} V,
\nabla^{\overline{\eta}}_W  F_{n+1}  ).
\end{eqnarray*}

\noindent Then, the equation above for $\overline{\eta} (V, \  R_{\overline{\eta}} ( W, F_{n+1} ) ( F_{n+1}) )$ becomes
\begin{eqnarray*}
&       &  \overline{\eta} (V, \  R_{\overline{\eta}} ( W, F_{n+1}
) ( F_{n+1}) )        \\   &      &      \\   &  =  &  -
\rho^{-1} ( \nabla^{\overline{g}}_W d \rho )(V) -  \overline{\eta}
( V,  {\rm II}^2 (W) ) +  F_{n+1}  \{  \overline{\eta} (V,  {\rm
II}(W) )  \}             \\   &      &      \\   &       &  + 2
\overline{\eta} ( \nabla^{\overline{\eta}}_{F_{n+1}} V,
\nabla^{\overline{\eta}}_W  F_{n+1}  ) +  \overline{\eta} ( {\rm
II}(W),  [ F_{n+1},  V ]  )       -  \overline{\eta} ( {\rm
II}(V),  [ F_{n+1},  W ]  )       \\   &      &       \\   & = &
F_{n+1}  \{ {\rm II} (V,  W) )  \}     -  \rho^{-1} (
\nabla^{\overline{g}}_W d \rho )(V) -  \overline{\eta} ( V,  {\rm
II}^2 (W) )       \\   &      &       \\   &     & - 2
\overline{\eta}  \Big(  \rho^{-1}  \nabla^{\eta}_{E_{n+1}} V +
\rho^{-1}  \Theta(V) -  {\rm II}(V) - \rho^{-1} [ \zeta,  V ],  \
{\rm II}(W)   \Big)          \\   &      &     \\   &      & +
\overline{\eta}  \Big(  \rho^{-1}  \nabla^{\eta}_{E_{n+1}} V +
\rho^{-1}  \Theta(V) - \rho^{-1} [ \zeta,  V ],  \   {\rm II}(W)
\Big)          \\   &      &     \\   &       &  -  \overline{\eta}  \Big(
\rho^{-1}  \nabla^{\eta}_{E_{n+1}} W + \rho^{-1}  \Theta(W) -
\rho^{-1} [ \zeta,  W ],  \   {\rm II}(V)   \Big)      .
\end{eqnarray*}

\noindent Rewriting yields,
\begin{eqnarray*}
&     &   E_{n+1} \{  {\rm II}(V, W)  \}  =   \rho  F_{n+1} \{
{\rm II}(V, W)  \}      +  \zeta \{  {\rm II}(V, W)  \}     \\   &
&       \\   &  =  &   \rho \cdot   \overline{\eta} (V, \
R_{\overline{\eta}} ( W, F_{n+1} ) ( F_{n+1}) )
 +   ( \nabla^{\overline{g}}_W d \rho )(V) - \rho  \overline{\eta} ( V,  {\rm II}^2 (W) )       \\
&      &       \\   &      &  - \overline{\eta} (  {\rm II}(V),  [
\zeta, W ] ) - \overline{\eta} (  {\rm II}(W),  [ \zeta, V ] )
+  \overline{\eta} (  {\rm II}(V),    \nabla^{\eta}_{E_{n+1}} W
)    +  \overline{\eta} (  {\rm II}(W),    \nabla^{\eta}_{E_{n+1}}
V   )               \\   &       &       \\   &      &
 +  \overline{\eta} (  {\rm II}(V),   \Theta(W) )  +     \overline{\eta} (  {\rm II}(W),   \Theta(V) )   +  \zeta \{  {\rm II}(V, W)  \}  .
\end{eqnarray*}

\noindent With the help of the equation (ii) in Lemma 3.1 and the identities,
\begin{eqnarray*}
E_{n+1} \{  {\rm
II}(V, W)  \} & = & ( \nabla^{\eta}_{E_{n+1}} {\rm II} ) (V, W)  +
{\rm II} (  \nabla^{\eta}_{E_{n+1}} V, W )   + {\rm II} (  V,
\nabla^{\eta}_{E_{n+1}} W ) ,       \\    &       &
\\      \zeta \{  {\rm II}(V, W)  \}  & =  &  (
\nabla^{\overline{g}}_{\zeta} {\rm II} ) (V, W)  + {\rm II} (
\nabla^{\overline{g}}_{\zeta} V, W )   + {\rm II} (  V,
\nabla^{\overline{g}}_{\zeta} W ) ,
\end{eqnarray*}
we obtain the asserted formula of the proposition.        $\rule{2mm}{2mm}$

\bigskip \noindent
{\bf Remark 3.1} Contracting the equation in Proposition 3.2 and applying (iii) in Lemma 3.1,
we obtain the following formula for ${\rm Ric}_{\overline{\eta}} (F_{n+1}, F_{n+1})$:
\begin{eqnarray*}
&     &  \rho \cdot {\rm Ric}_{\overline{\eta}} (F_{n+1}, F_{n+1})
  =  {\rm Tr}_{\overline{g}} ( \nabla^{\eta}_{E_{n+1}} {\rm II})  - {\rm Tr}_{\overline{g}} ( \nabla^{\overline{g}}_{\zeta}
{\rm II} )    + \rho \cdot {\rm Tr}_{\overline{g}} ( {\rm II}^2 )  \\
&     &      \\                                                                                &     & \qquad \qquad - 2 \sum_{i=1}^n  \overline{\eta} ( {\rm II} (F_i),  \,  \nabla^{\overline{g}}_{F_i} \zeta )   - 2 \sum_{i=1}^n \overline{\eta} ( \Theta(F_i), \, {\rm II}(F_i) )- \sum_{i=1}^n ( \nabla^{\overline{g}}_{F_i} d \rho ) (F_i).
\end{eqnarray*}

\bigskip \noindent  \section{Representation of the Dirac equation on manifolds with codimension one foliation }

\noindent
In this section we will represent the Dirac equation on $(Q^{n+1}, \overline{\eta})$
hyperbolically with respect to codimension one foliation (see
Proposition 4.1 and Corollary 4.1).
Let us fix a slice $(M^n,  \overline{g}
)$ of the foliated manifold $(
Q^{n+1},  \overline{\eta} )$.
We will identify $\Sigma(Q)_{\overline{\eta}}$ with $\Sigma(Q)_{\eta}$ and $\psi \in \Gamma (\Sigma(Q)_{\overline{\eta}} )$ with
its pullback $\widetilde{K}(\psi)$, via the natural isomorphism
$ \widetilde{K} : \Sigma(Q)_{\overline{\eta}} \longrightarrow \Sigma(Q)_{\eta} $, and write
simply as $\Sigma(Q)$ and $\psi$, respectively.                                                              Depending on the
dimension $n+1$ of the manifold $Q^{n+1}$, we will use two different
Clifford multiplications in the subbundle $\Sigma (M) \subset \Sigma(Q)$. For the realization
of the Clifford algebra over ${\mathbb R}$, we refer to [11].  Let ${\rm Cl}(M)$ (resp.
${\rm Cl}(Q)$) denote the Clifford bundle over $M^n$ (resp. $Q^{n+1}$).   \par (i)
In case of $n = 2m$, we use the Clifford multiplication  $ {\rm
Cl}(M)  \times  \Sigma(M)   \longrightarrow \Sigma(M) $ that is
naturally related to the one  $ {\rm Cl}(Q)  \times  \Sigma(Q)
\longrightarrow \Sigma(Q) $  via
\begin{eqnarray*}
\pi_{\ast} (F_i  \cdot  \psi)  & =  & F_i   \cdot  ( \pi_{\ast}
\psi ) ,   \quad  1  \leq i  \leq  2m ,         \\ &      &
\\    \pi_{\ast} (F_{2m+1}  \cdot  \psi)  & =  &  ( \sqrt{-1}
)^{m+1}   \mu_{\overline{g}}   \cdot   ( \pi_{\ast} \psi ) ,
\end{eqnarray*}
where $\pi_{\ast} :  \Sigma(Q)   \longrightarrow   \Sigma(M)$ is
the restriction map and $\mu_{\overline{g}}$ is the volume element of $(M^{2m}, \overline{g})$.
The second relation is an immediate consequence of the algebraic relation
\[      F_{2m+1} \cdot \psi = ( \sqrt{-1} )^{m+1} F_1  \cdots  F_{2m} \cdot \psi .      \]
\indent (ii) In the other case  $n = 2m -1$,
we identify the spinor bundle $\Sigma(M)$ with the positive part
$\Sigma^+ (Q)$ of the bundle $\Sigma (Q)$ restricted to $M^{2m-1}$,
and we use the Clifford multiplication $ {\rm Cl}(M)  \times
\Sigma(M)   \longrightarrow \Sigma(M) $ that is naturally related
to $ {\rm Cl}^+(Q)  \times  \Sigma^+(Q)   \longrightarrow
\Sigma^+(Q) $  via
\begin{eqnarray*}
\pi_{\ast}^+ (F_i  \cdot F_{2m} \cdot  \psi^+)  & =  &    F_i
\cdot  ( \pi_{\ast}^+ \psi^+ ) ,   \quad  1  \leq i  \leq  2m-1,
\\   &      &     \\     \pi_{\ast}^+ (F_k  \cdot  F_l  \cdot  \psi^+)  &
=  &  F_k  \cdot  F_l  \cdot   ( \pi_{\ast}^+ \psi^+ ) ,     \quad   1
\leq k < l   \leq 2m-1 ,
\end{eqnarray*}
where $\pi_{\ast}^+ :  \Sigma^+(Q)   \longrightarrow   \Sigma(M)$ is
the restriction map and ${\rm Cl}^+(Q)$ is the positive part of ${\rm Cl}(Q)$.

\bigskip
Recall that the spin derivatives $\nabla^{\overline{\eta}} \psi, \, \nabla^{\overline{g}} \psi$
are related, on $Q^{n+1}$, by
\[
\nabla^{\overline{\eta}}_V \psi   =  \nabla^{\overline{g}}_V
\psi  +  \frac{1}{2}  {\rm II}(V)  \cdot  F_{n+1}
\cdot  \psi .      \]

\bigskip  \noindent
In view of  the rule of Clifford multiplication described above, we find that,
 in case of $n = 2m$, the formula is projected to the slice $(M^{2m}, \overline{g})$ as
\[          \pi_{\ast} ( \nabla^{\overline{\eta}}_V \psi )   =
\nabla^{\overline{g}}_V  ( \pi_{\ast} \psi )  +  \frac{1}{2}
( \sqrt{-1} )^{m+1}   {\rm II}(V)  \cdot
\mu_{\overline{g}}  \cdot  ( \pi_{\ast}  \psi  ) .       \]              However, in
the other case $n = 2m-1$, the projection is only possible if  $\psi = \psi^+$
belongs entirely to the positive part $\Sigma^+(Q) $ of $\Sigma(Q)$,
the projected formula being given by    \[                  \pi_{\ast}^+ (
\nabla^{\overline{\eta}}_V \psi^+ )  =  \nabla^{\overline{g}}_V (
\pi_{\ast}^+  \psi^+ )  +  \frac{1}{2}  {\rm II}(V)
\cdot  ( \pi_{\ast}^+ \psi^+ ) .        \]             Nevertheless,
we may regard not only $\nabla^{\overline{g}}_V  \psi^+$ but also   $\nabla^{\overline{g}}_V  \psi^- , \, \psi^- \in  \Gamma( \Sigma^-(Q) )$ (e.g. $\psi^- =
F_{2m}  \cdot  \psi^+$), as well-defined spinor fields on $Q^{2m}$, not
projected to the slice $M^{2m-1}$. Therefore, the following formula
makes sense.

\bigskip  \noindent
\begin{lem}
\[       \nabla^{\overline{g}}_V   (  F_{2m}
\cdot  \psi^+   )   =   F_{2m}  \cdot  \nabla^{\overline{g}}_V
\psi^+   .        \]
\end{lem}

\noindent
{\bf Proof.}  We compute
\begin{eqnarray*}
&       &  \nabla^{\overline{\eta}}_V   (  F_{2m}  \cdot  \psi^+
)    =   \nabla^{\overline{\eta}}_V    F_{2m}  \cdot  \psi^+ +
F_{2m}  \cdot   \nabla^{\overline{\eta}}_V   \psi^+         \\  &
&          \\   & =   &  -  {\rm II} (V)  \cdot  \psi^+   +  F_{2m}
\cdot \Big\{ \nabla^{\overline{g}}_V  \psi^+  +  \frac{1}{2}
  {\rm II}(V)  \cdot  F_{2m}  \cdot  \psi^+   \Big\}
\\   &       &       \\   & = &  F_{2m}  \cdot
\nabla^{\overline{g}}_V   \psi^+ -  \frac{1}{2}   {\rm II} (V)
\cdot  \psi^+ .
\end{eqnarray*}
On the other hand,
\begin{eqnarray*}
&       &  \nabla^{\overline{\eta}}_V   (  F_{2m}  \cdot  \psi^+
)   =
 \nabla^{\overline{g}}_V   (  F_{2m}  \cdot  \psi^+   ) +  \frac{1}{2}    {\rm II} (V)  \cdot  F_{2m}  \cdot  F_{2m}
\cdot  \psi^+        \\   &      &       \\    &   =  &
\nabla^{\overline{g}}_V   (  F_{2m}  \cdot  \psi^+   )  -
\frac{1}{2}    {\rm II} (V)   \cdot  \psi^+    .
\end{eqnarray*}
Comparing the latter equation with the former, we complete the proof.
      $\rule{2mm}{2mm}$

\bigskip
In order to represent the Dirac equation
\begin{eqnarray*}
&      & D_{\overline{\eta}} \psi  =   \lambda_Q
\psi              \\   &       &     \\
& =   & \sum_{i=1}^n  F_i \cdot  \nabla^{\overline{g}}_{F_i} \psi
- \frac{1}{2}  ( {\rm Tr}_{\overline{g}} {\rm II} )
F_{n+1} \cdot \psi  +   F_{n +1} \cdot
\nabla^{\overline{\eta}}_{F_{n+1}} \psi  ,
\end{eqnarray*}
with respect to reference metric, we need the following lemma that one verifies
straightforwardly using Proposition 3.1.

\bigskip \noindent
\begin{lem}
\begin{eqnarray*}
&     & \rho  \cdot \Big\{  \widetilde{K} \circ
\nabla^{\overline{\eta}}_{F_{n+1}}  \circ  \big( \widetilde{K}
\big)^{-1}   \Big\} (\varphi)      \\   &     &     \\   &  =  &
\nabla^{\eta}_{E_{n+1}} \varphi - \nabla^g_{\zeta} \varphi   -
\frac{\rho}{4}  \sum_{i=1}^n  E_i \cdot (L \circ {\rm II} \circ
L^{-1})(E_i) \cdot \varphi         \\   &     &     \\   &      & -
\frac{1}{4} \sum_{i=1}^n E_i \cdot ( \nabla^{\eta}_{E_{n+1}} L )
(L^{-1} E_i) \cdot \varphi + \frac{1}{4} \sum_{i=1}^n E_i \cdot (
\nabla^g_{\zeta} L ) (L^{-1} E_i) \cdot \varphi          \\   &
&     \\   &    & +   \frac{1}{4} \sum_{i=1}^n E_i \cdot L(
\nabla^g_{L^{-1}E_i} \zeta ) \cdot \varphi
 +   \frac{1}{4} \sum_{i=1}^n  E_i \cdot (L \circ \Theta \circ L^{-1})(E_i) \cdot \varphi          \\
&      &     \\   &      & + \frac{1}{2} \sum_{i=1}^n d \rho
(L^{-1} E_i) E_i \cdot E_{n+1} \cdot \varphi.
\end{eqnarray*}
\end{lem}

\bigskip \noindent
Lemma 4.2, combined with Proposition 2.3, yields the following hyperbolic representation of the Dirac equation immediately.

\bigskip  \noindent
\begin{pro}
\begin{eqnarray*}
  \nabla^{\eta}_{E_{n+1}} \psi  &  = &    \nabla^g_{\zeta} \psi  - \lambda_Q \rho
 E_{n+1} \cdot \psi     +  \frac{\rho}{2} (
{\rm Tr}_{\overline{g}} {\rm II} ) \psi          \\   &      &
\\   &       &  + \rho E_{n+1}  \cdot  \Big\{  \sum_{i=1}^n  E_i
\cdot   \nabla^g_{L^{-1} E_i} \psi +  \frac{1}{4}
\sum_{j,k,l=1}^n ( \Lambda_g )_{jkl} E_j  \cdot E_k  \cdot E_l
\cdot  \psi    \Big\}      \\   &       &       \\   &      &   +
\frac{1}{4} \sum_{i=1}^n E_i \cdot ( \nabla^{\eta}_{E_{n+1}} L )
(L^{-1} E_i) \cdot \psi + \frac{\rho}{4}   \sum_{i=1}^n  E_i
\cdot (L \circ {\rm II} \circ L^{-1})(E_i) \cdot \psi     \\   &
&     \\   &    &  - \frac{1}{4} \sum_{i=1}^n E_i \cdot (
\nabla^g_{\zeta} L ) (L^{-1} E_i) \cdot \psi      -
\frac{1}{4} \sum_{i=1}^n E_i \cdot L( \nabla^g_{L^{-1}E_i} \zeta )
\cdot \psi
   \\
&      &     \\   &      &   -   \frac{1}{4} \sum_{i=1}^n  E_i
\cdot (L \circ \Theta \circ L^{-1})(E_i) \cdot \psi    -
\frac{1}{2} \sum_{i=1}^n d \rho (L^{-1} E_i) E_i \cdot E_{n+1}
\cdot \psi.
\end{eqnarray*}
\end{pro}

\bigskip \noindent
Although the equation in Proposition 4.1 is valid in both cases, $n = 2m$ and $n=2m-1$,
 it is also very useful, in the latter case $n=2m-1$, to consider the decomposition of spinor fields,
\[           \psi = \psi^+ + F_{2m} \cdot \varphi^+,  \qquad \psi^+, \varphi^+  \in \Gamma( \Sigma^+(Q) ) ,         \]
and rewrite the representation in Proposition 4.1 equivalently as follows.

\bigskip \noindent
\begin{cor}
\begin{eqnarray*}
\nabla^{\eta}_{E_{2m}} \psi^+
 &  = &    \nabla^g_{\zeta} \psi^+  + \lambda_Q \rho
 \varphi^+     +  \frac{\rho}{2} (
{\rm Tr}_{\overline{g}} {\rm II} ) \psi^+          \\   &      &
\\   &       &  + \rho E_{2m}  \cdot  \Big\{  \sum_{i=1}^{2m-1}  E_i
\cdot   \nabla^g_{L^{-1} E_i} \psi^+  +  \frac{1}{4}
\sum_{j,k,l=1}^{2m-1} ( \Lambda_g )_{jkl} E_j  \cdot E_k  \cdot E_l
\cdot  \psi^+    \Big\}      \\   &       &       \\   &      &   +
\frac{1}{4} \sum_{i=1}^{2m-1} E_i \cdot ( \nabla^{\eta}_{E_{2m}} L )
(L^{-1} E_i) \cdot \psi^+ + \frac{\rho}{4}   \sum_{i=1}^{2m-1}  E_i
\cdot (L \circ {\rm II} \circ L^{-1})(E_i) \cdot \psi^+     \\   &
&     \\   &    &  - \frac{1}{4} \sum_{i=1}^{2m-1} E_i \cdot (
\nabla^g_{\zeta} L ) (L^{-1} E_i) \cdot \psi^+      -
\frac{1}{4} \sum_{i=1}^{2m-1} E_i \cdot L( \nabla^g_{L^{-1}E_i} \zeta )
\cdot \psi^+
   \\
&      &     \\   &      &   -   \frac{1}{4} \sum_{i=1}^{2m-1}  E_i
\cdot (L \circ \Theta \circ L^{-1})(E_i) \cdot \psi^+    -
\frac{1}{2} \sum_{i=1}^{2m-1} d \rho (L^{-1} E_i) E_i \cdot E_{2m}
\cdot \psi^+,        \\
&       &     \\
&       &     \\
  \nabla^{\eta}_{E_{2m}} \varphi^+    &  = &    \nabla^g_{\zeta} \varphi^+  - \lambda_Q \rho
 \psi^+     +  \frac{\rho}{2} (
{\rm Tr}_{\overline{g}} {\rm II} ) \varphi^+          \\   &      &
\\   &       &  - \rho E_{2m}  \cdot  \Big\{  \sum_{i=1}^{2m-1}  E_i
\cdot   \nabla^g_{L^{-1} E_i} \varphi^+ +  \frac{1}{4}
\sum_{j,k,l=1}^{2m-1} ( \Lambda_g )_{jkl} E_j  \cdot E_k  \cdot E_l
\cdot  \varphi^+    \Big\}      \\   &       &       \\   &      &   +
\frac{1}{4} \sum_{i=1}^{2m-1} E_i \cdot ( \nabla^{\eta}_{E_{2m}} L )
(L^{-1} E_i) \cdot \varphi^+ + \frac{\rho}{4}   \sum_{i=1}^{2m-1}  E_i
\cdot (L \circ {\rm II} \circ L^{-1})(E_i) \cdot \varphi^+     \\   &
&     \\   &    &  - \frac{1}{4} \sum_{i=1}^{2m-1} E_i \cdot (
\nabla^g_{\zeta} L ) (L^{-1} E_i) \cdot \varphi^+      -
\frac{1}{4} \sum_{i=1}^{2m-1} E_i \cdot L( \nabla^g_{L^{-1}E_i} \zeta )
\cdot \varphi^+
   \\
&      &     \\  &      &   -   \frac{1}{4} \sum_{i=1}^{2m-1}  E_i
\cdot (L \circ \Theta \circ L^{-1})(E_i) \cdot \varphi^+    +
\frac{1}{2} \sum_{i=1}^{2m-1} d \rho (L^{-1} E_i) E_i \cdot E_{2m}
\cdot \varphi^+.
\end{eqnarray*}
\end{cor}

\bigskip
We close this section with representing the energy-momentum tensor
\[         T_{\overline{\eta}} (X, Y) = \frac{\epsilon}{4} \Big( X \cdot \nabla^{\overline{\eta}}_Y
\psi +  Y \cdot \nabla^{\overline{\eta}}_X \psi, \ \psi \Big),   \qquad \epsilon = \pm 1 .
\]
hyperbolically with respect to codimension one foliation.
To this end, it is important to notice that, if
$\psi$ be a
solution of the Dirac equation $ D_{\overline{\eta}} \psi =
\lambda_Q \psi, \ \lambda_Q \in {\mathbb R}$, then the following equation is valid :
\[
\nabla^{\overline{\eta}}_{F_{n+1}} \psi   =  - \lambda_Q   F_{n+1}  \cdot  \psi  +  \frac{1}{2} ( {\rm
Tr}_{\overline{g}} {\rm II} )  \psi  + F_{n+1}  \cdot (
\sum_{i=1}^n F_i  \cdot  \nabla^{\overline{g}}_{F_i} \psi ) .
\]

\bigskip  \noindent
\begin{pro}       For any solution $\psi$ of the Dirac
equation  $D_{\overline{\eta}} \psi  =  \lambda_Q  \psi $ on $(Q^{n+1}, \overline{\eta})$, we have
\begin{eqnarray*}
{\rm Tr}_{\overline{\eta}} ( T_{\overline{\eta}} )  & = &  \frac{\epsilon
\lambda_Q}{2} ( \psi, \  \psi )  ,   \\   &        &
\cr T_{\overline{\eta}} (V, W) & = &  \frac{\epsilon}{4} \Big(
V \cdot \nabla^{\overline{g}}_W \psi + W \cdot
\nabla^{\overline{g}}_V \psi  ,  \ \psi  \Big)        \\
&       &         \\  &       &      +   \frac{\epsilon}{8} \Big(
\{ V \cdot  {\rm II}(W)  + W \cdot  {\rm II}(V)  \} \cdot  F_{n+1}
\cdot \psi ,    \  \psi    \Big)     ,      \\  &       &
\cr T_{\overline{\eta}} (V, F_{n+1}) & = & \frac{\epsilon}{4} \Big(
F_{n+1} \cdot  \Big\{ \nabla^{\overline{g}}_V  \psi - V \cdot
 ( \sum_{i=1}^n F_i  \cdot  \nabla^{\overline{g}}_{F_i} \psi )  \Big\} ,  \   \psi  \Big),       \\
&        &       \\  T_{\overline{\eta}} (F_{n+1}, F_{n+1}) & = & -
\frac{\epsilon}{2} \Big( ( \sum_{i=1}^n F_i  \cdot
\nabla^{\overline{g}}_{F_i} \psi )  -  \lambda_Q   \psi  ,  \  \psi
\Big).
\end{eqnarray*}
\end{pro}

\bigskip  \noindent
In case of $n=2m-1$, we consider the decomposition $\psi = \psi^+ + F_{2m} \cdot \varphi^+$
and can equivalently rewrite the formulas in Proposition 4.2 as follows.

\bigskip \noindent
\begin{cor}                   For any solution $\psi$ of the Dirac
equation  $D_{\overline{\eta}} \psi  =   \lambda_Q
\psi $ on $(Q^{2m},  \overline{\eta} )$
, where $\psi = \psi^+ + F_{2m} \cdot \varphi^+$, we have
\begin{eqnarray*}
{\rm Tr}_{\overline{\eta}} ( T_{\overline{\eta}} )  & = & \frac{\epsilon
\lambda_Q}{2}  \Big\{ ( \psi^+, \psi^+ ) + ( \varphi^+,  \varphi^+ ) \Big\}  ,   \\                                     &   &        \\                                                                      T_{\overline{\eta}} (V, W) & = & \frac{\epsilon}{4} \Big(
 V \cdot \nabla^{\overline{g}}_W \psi^+ + W \cdot
\nabla^{\overline{g}}_V \psi^+  ,  \ F_{2m} \cdot \varphi^+  \Big)
\\                                                                                             &       &         \\
&       &  + \frac{\epsilon}{4} \Big(
V \cdot \nabla^{\overline{g}}_W \varphi^+ + W \cdot
\nabla^{\overline{g}}_V \varphi^+ ,  \ F_{2m} \cdot \psi^+ \Big)
\\                                                                                             &       &         \\                                                                            &       &      +   \frac{\epsilon}{4}
\Big( \{ V \cdot  {\rm II}(W)  + W \cdot  {\rm
II}(V)  \} \cdot \psi^+ ,    \  \varphi^+    \Big)         \\
&       &       \\
&       &      +   \frac{\epsilon}{2}
{\rm II}(V, W) ( \psi^+, \ \varphi^+  )    ,     \\
&       &       \\                                                                                                                                        T_{\overline{\eta}} (V, F_{2m}) & = &
\frac{\epsilon}{4} \Big(
\nabla^{\overline{g}}_V  \psi^+ - V \cdot
 ( \sum_{i=1}^{2m-1} F_i  \cdot  \nabla^{\overline{g}}_{F_i} \psi^+ )   ,  \   \varphi^+  \Big)       \\
&       &       \\                                                                                                                                         &       & -
\frac{\epsilon}{4} \Big(
\nabla^{\overline{g}}_V  \varphi^+ - V \cdot
 ( \sum_{i=1}^{2m-1} F_i  \cdot  \nabla^{\overline{g}}_{F_i} \varphi^+ )   ,  \   \psi^+  \Big) ,       \\
&        &       \cr T_{\overline{\eta}} (F_{2m}, F_{2m}) & = &
\frac{\epsilon}{2} \Big( F_{2m} \cdot  (
\sum_{i=1}^{2m-1} F_i  \cdot  \nabla^{\overline{g}}_{F_i} \psi^+ )
  ,  \  \varphi^+  \Big)    \\
&     &     \\
&     &  + \frac{\epsilon}{2} \Big( F_{2m} \cdot  (
\sum_{i=1}^{2m-1} F_i  \cdot  \nabla^{\overline{g}}_{F_i} \varphi^+ )
  ,  \  \psi^+  \Big)   \\
&     &     \\
&     &   + \frac{\epsilon
\lambda_Q}{2}  \Big\{ ( \psi^+, \psi^+ ) + ( \varphi^+,  \varphi^+ ) \Big\}  .  \end{eqnarray*}   \end{cor}

\bigskip \noindent
\section{A sufficient condition for the existence of solutions to the weak Killing equation}

\noindent
Let us suppose that $(Q^{n+1}, \overline{\eta})$ satisfies
\[
{\rm Ric}_{\overline{\eta}} (V, W) = \frac{S_{\overline{\eta}}}{2} \overline{\eta}(V,W),
\quad {\rm Ric}_{\overline{\eta}} (V, F_{n+1}) = 0,
\quad d S_{\overline{\eta}} (V) = 0
\]
for all $V, W \in E_{n+1}^{\bot}$. Then the weak Killing equation becomes
\[
\nabla^{\overline{\eta}}_V \psi  =  \nabla^{\overline{g}}_V \psi + \frac{1}{2} {\rm II}(V)
\cdot F_{n+1} \cdot \psi   = \frac{d S_{\overline{\eta}} (F_{n+1})}{2n S_{\overline{\eta}}}
V \cdot F_{n+1} \cdot \psi
\]
and
\begin{eqnarray*}
\nabla^{\overline{\eta}}_{F_{n+1}} \psi & = & - \lambda_Q F_{n+1} \cdot \psi
+ \frac{(n+1) d S_{\overline{\eta}}(F_{n+1}) }{ 2n S_{\overline{\eta}} } \psi  +
\frac{d S_{\overline{\eta}}(F_{n+1}) }{ 2n S_{\overline{\eta}} }  F_{n+1} \cdot F_{n+1} \cdot \psi      \\
&      &      \\
&  =  & - \lambda_Q F_{n+1} \cdot \psi + \frac{d S_{\overline{\eta}}(F_{n+1}) }{ 2 S_{\overline{\eta}} } \psi .
\end{eqnarray*}

\bigskip \noindent
Thus we have proved the following proposition.

\bigskip \noindent
\begin{pro}    Let $(Q^{n+1}, \overline{\eta} )$ satisfy the following conditions:
\begin{eqnarray*}
&     &   {\rm Ric}_{\overline{\eta}} (V, W) = \frac{S_{\overline{\eta}}}{2} \overline{\eta}(V,W),
 \quad  {\rm Ric}_{\overline{\eta}} (V, F_{n+1}) = 0,
\quad d S_{\overline{\eta}} (V) = 0   ,    \\
&     &       \\
&     &   {\rm II} (V, W) = \frac{d S_{\overline{\eta}}(F_{n+1}) }{ n S_{\overline{\eta}} }
\overline{\eta}(V, W).
\end{eqnarray*}
Under this assumption the weak Killing equation is equivalent to the system of  differential equations,
\[
\nabla^{\overline{g}}_V \psi  =  0   \qquad   \mbox{and}  \qquad
\nabla^{\overline{\eta}}_{F_{n+1}} \psi
 = - \lambda_Q F_{n+1} \cdot \psi + \frac{1}{2} {\rm Tr}_{\overline{g}}( {\rm II} ) \psi .
\]
\end{pro}

\noindent
As an application of Proposition 5.1, we are going to prove below that every parallel spinor
may evolve to a WK-spinor (Theorem 5.1). For this purpose, we first show that there indeed exist
some special metrics satisfying the hypothesis of Proposition 5.1.
Let $Q^{n+1} = M^n \times {\mathbb R}$ be a product manifold, and let the product metric $\eta = g_M \times g_{\mathbb R}$
be the reference metric on $Q^{n+1}$, where  $g_M$ indicates an arbitrary Riemannian metric
on $M^n$ and $g_{\mathbb R}$ the standard metric on the real line ${\mathbb R}$.
We write $g_{\mathbb R} = dt \otimes dt$, using the standard coordinate $t \in {\mathbb R}$.                                                  By $(E_1, \ldots, E_n)$ we denote a local orthonormal frame on $(M^n, g_M)$ as well as its lift to $(Q^{n+1}, \eta)$.   Let $E_{n+1} = \frac{d}{dt}$ denote the unit vector field on $({\mathbb R}, g_{\mathbb R})$ as well as the lift to $(Q^{n+1}, \eta)$.
Then it is clear that
$ \Theta(V) = - \nabla^{\eta}_V E_{n+1} = 0  $
for all vector fields $V \in  E_{n+1}^{\bot}$.
For simplicity, we denote by ${\rm WP}(g_M; a)$ the following class of metrics
(the {\it warped products of $g_M$ and $g_{\mathbb R}$ })  :                                             \[   \overline{\eta} = e^f \Big( \sum_{i=1}^n E^i \otimes E^i \Big) +  e^{af} dt \otimes
dt ,  \]
where $f : {\mathbb R}  \longrightarrow {\mathbb R}$ is a real-valued function and $a \in {\mathbb R}$
is a real number.

\bigskip  \noindent
\begin{lem}                 For  all $\overline{\eta} \in {\rm WP}( g_M ; a)$, we have:
\begin{eqnarray*}
{\rm II} (V, W) & = & - \frac{e^{-(\frac{a}{2} -1)f} f_t}{2}  \eta(V, W) ,   \\
&      &      \\
\overline{\eta} (V, {\rm II}^2 (W) ) & = &  \frac{e^{-(a-1)f} f_t f_t}{4}
\eta(V, W)   ,    \\
&     &     \\
{\rm Tr}_{\overline{g}} ({\rm II}) & = & - \frac{n e^{-\frac{a}{2}f} f_t}{2 }  ,  \\
&     &     \\
{\rm Tr}_{\overline{g}} ({\rm II}^2) & = & \frac{n e^{-af} f_t f_t}{4}  ,   \\
&     &     \\
( \nabla^{\eta}_{E_{n+1}} {\rm II} ) ( V, W )  & = & \Big\{ -  \frac{e^{-(\frac{a}{2} -1)f} f_{tt}}{2}
 + \frac{(a-2) e^{-(\frac{a}{2} -1)f} f_t f_t}{4} \Big\}
 \eta(V,W),
\end{eqnarray*}  where we have used the shorthand notation $f_t : = df(E_{n+1})$ and $\displaystyle f_{tt} : = \frac{d^2 f}{dt^2} = ( \nabla^{\eta}_{E_{n+1}} df) (E_{n+1}) $.
\end{lem}

\noindent
{\bf Proof.}  Since $L = e^{\frac{f}{2}} I$ (in the notations of Section 3), we have
\[         ( \nabla^g_V L)(W) = 0   \qquad  \mbox{and}  \qquad
( \nabla^{\eta}_{E_{n+1}} L)(W) = \frac{e^{\frac{f}{2}} f_t}{2} W ,             \]                and hence Corollary 3.1 gives
\[       {\rm II} (V, W)  = - \frac{e^{-(\frac{a}{2} -1)f} f_t}{2}  \eta(V, W) .        \]
Using this, one checks all the equations of the lemma easily.    
$\rule{2mm}{2mm}$

\bigskip \noindent
Substituting Lemma 5.1 in Proposition 3.2 as well as in the Gauss-Codazzi equations (Lemma 3.1),
we obtain the following lemma immediately.

\bigskip  \noindent
\begin{lem}                 For all $\overline{\eta} \in {\rm WP}( g_M ; a)$, we have:
\begin{eqnarray*}
{\rm Ric}_{\overline{\eta}} (V, W) & = &  {\rm Ric}_g (V, W) + \Big\{ -  \frac{f_{tt}}{2e^{(a-1)f}}
+  \frac{(a-n) f_t f_t}{4 e^{(a-1)f}} \Big\} \eta(V, W) ,      \\
&     &    \\
{\rm Ric}_{\overline{\eta}} (V, F_{n+1} ) & = & 0,     \\
&     &     \\
{\rm Ric}_{\overline{\eta}}  ( F_{n+1}, F_{n+1} ) & = &
- \frac{n f_{tt}}{2 e^{af}}  +  \frac{n(a-1) f_t f_t}{4 e^{af}}
 ,        \\
&     &     \\
S_{\overline{\eta}} & = &  e^{-f} S_g - \frac{n f_{tt}}{e^{af}}
+ \frac{n (2a-n-1) f_t f_t}{4 e^{af}}  .
\end{eqnarray*}
\end{lem}

\bigskip \noindent
\begin{lem}   \   $\overline{\eta} \in {\rm WP}( g_M ; a)$ satisfies the
hypothesis of Proposition 5.1
if and only if
either   \[    f (t) = ct  \qquad \Big( a = \frac{n}{2} \Big)       \]   or
 \[     f(t) =  \frac{4}{n-2a}    \log \Big\{ 1 + \frac{(n-2a) ct}{4} \Big\}    \qquad                     \Big( a \neq  \frac{n}{2},  \ 1 + \frac{(n-2a) ct}{4} > 0   \Big) ,       \]
where $f(0) = 0$ and $c := f_t(0)$.   The scalar curvature of the metric $\overline{\eta}$ is non-positive
and is equal to
\[         - \frac{n c^2}{4} e^{- \frac{nc}{2}t}    \qquad   \mbox{and}  \qquad
- \frac{nc^2}{4}  \Big\{ \frac{4 + (n-2a) ct}{4}  \Big\}^{ \frac{-2n}{n-2a} } ,
\]           respectively.
\end{lem}

\noindent
{\bf Proof.}  Using Lemma 5.1 and 5.2, we compute
\[
{\rm Ric}_{\overline{\eta}} (V, W) - \frac{1}{2} S_{\overline{\eta}} \overline{\eta} (V, W)
 =  \frac{(n-1)}{8} e^{-af} \Big\{  4 f_{tt} + (n-2a) f_t f_t  \Big\} \overline{\eta}(V, W)
\]
and
\begin{eqnarray*}
&     &    n S_{\overline{\eta}}  {\rm II} (V, W)- d S_{\overline{\eta}} (F_{n+1})
\overline{\eta}(V, W)      \\
&     &      \\
&  =  & \frac{n}{8}  e^{- \frac{3}{2} af} \Big\{ 8 f_{ttt} - 4(4a-2n-1) f_t f_{tt}
+ (2a-n)(2a-n-1) f_t f_t f_t   \Big\} \overline{\eta} (V, W) .
\end{eqnarray*}
Since $4 f_{tt} + (n-2a) f_t f_t = 0$ implies
\[      8 f_{ttt} - 4 (4a-2n-1) f_t f_{tt} + (2a-n)(2a-n-1) f_t f_t f_t = 0 ,       \]
the hypothesis of Proposition 5.1 is satisfied if and only if
$ 4 f_{tt} + (n-2a) f_t f_t = 0 ,$
which can be solved completely as given in the lemma.   $\rule{2mm}{2mm}$

\bigskip \noindent
\begin{pro}        Let $(M^{2m}, g_M)$ be a Riemannian manifold admitting
 parallel spinors. Then, for any real number $\lambda_Q \in {\mathbb R}$, there exists
a warped product metric $\overline{\eta} \in {\rm WP} (g_M;a)$ (see Lemma 5.3) on $Q^{2m+1} = M^{2m} \times
{\mathbb R}$ such that $( Q^{2m+1}, \overline{\eta} )$ admits a WK-spinor to WK-number
$\lambda_Q$.
\end{pro}

\noindent
{\bf Proof.} Let $\psi_M^+ \in \Gamma ( \Sigma^+(M) )$ be a parallel spinor on $M^{2m}$.
Let $ \psi^+ = h^+ \psi_M^+$ be a spinor field on $Q^{2m+1}$ defined by
\[     \psi^+(x, t) = h^+(t) \psi_M^+(x),  \qquad (x , t) \in M^{2m} \times {\mathbb R} ,
\]
where $h^+ : {\mathbb R} \longrightarrow {\mathbb C}$ is a complex-valued function with
$h^+(0) = 1$.
Let $\overline{\eta} \in {\rm WP} (g_M;a)$ be a warped product metric given as in Lemma 5.3.
Now recall that,
 in our realization of Clifford algebra, the volume form $\mu_M$ of $(M^{2m}, g_M)$ acts on $\Sigma^{\pm}(M)$
via    \[               \mu_M \cdot \psi^+_M =   ( \sqrt{-1} )^m  \psi^+_M  ,  \qquad
 \mu_M \cdot \psi^-_M =   - ( \sqrt{-1} )^m  \psi^-_M    .                   \]
Since $E_{2m+1} \cdot \varphi = ( \sqrt{-1} )^{m+1} \mu_M \cdot \varphi$
for all spinor fields $\varphi$ on $Q^{2m +1}$, we have
\[            E_{2m+1} \cdot \psi^+_M =  ( \sqrt{-1} )^{2m+1}  \psi^+_M ,   \qquad                        E_{2m+1} \cdot \psi^-_M =  - ( \sqrt{-1} )^{2m+1}  \psi^-_M  .            \]
By Proposition 5.1, $\psi^+$ is a WK-spinor to WK-number $\lambda_Q$ if and only if the function $h^+$ satisfies
\[          h^+_t = -( \sqrt{-1} )^{2m+1} \lambda_Q e^{\frac{a}{2} f} h^+ - \frac{m}{2} f_t h^+ ,     \]
which obviously allows a global solution.   $\rule{2mm}{2mm}$

\bigskip \noindent
We now extend Proposition 5.2 so as to include the case $n=2m-1$.

\bigskip \noindent
\begin{thm}         Let $(M^n, g_M)$ be a Riemannian manifold admitting
 parallel spinors. Then, for any real number $\lambda_Q \in {\mathbb R}$, there exists
a warped product metric $\overline{\eta} \in {\rm WP} (g_M;a)$ (see Lemma 5.3) on $Q^{n+1} = M^{n} \times
{\mathbb R}$ such that $( Q^{n+1}, \overline{\eta} )$ admits a WK-spinor to WK-number
$\lambda_Q$.
\end{thm}

\noindent
{\bf Proof.} Because of Proposition 5.2, it suffices to prove the theorem for the case $n=2m-1$.   Let $\varphi_M^+ \in \Gamma ( \Sigma^+(Q) )$ be a parallel spinor on $M^{2m-1}$, and let
$\varphi = h^+ \varphi^+_M + k^+ E_{2m} \cdot \varphi^+_M$ be a spinor field on $Q^{2m}$ defined by
\[      \varphi(x,t) = h^+(t) \varphi^+_M(x) + k^+(t) E_{2m} \cdot \varphi^+_M (x),  \qquad
(x,t) \in M^{2m-1} \times {\mathbb R},       \]
where $h^+ , k^+: {\mathbb R} \longrightarrow {\mathbb C}$ are complex-valued functions with
$h^+(0) = k^+(0) = 1$. Let $\overline{\eta} \in {\rm WP} (g_M;a)$ be a warped product metric given as in Lemma 5.3. By Proposition 5.1, $\varphi$ is a WK-spinor to WK-number $\lambda_Q$ if and only if
$(h^+, k^+)$ satisfies the system of differential equations,
\[
h^+_t  =   - \frac{2m-1}{4} f_t h^+ + \lambda_Q e^{\frac{a}{2}f} k^+    \quad   \mbox{and}  \quad
k^+_t  =  - \lambda_Q e^{\frac{a}{2}f} h^+ - \frac{2m-1}{4} f_t k^+ .
\]
This is a linear homogeneous system and hence allows a global solution.   
$\rule{2mm}{2mm}$

\bigskip \noindent
{\bf Remark 5.1} The WK-spinors constructed on ${\mathbb R}^3$ at the end of Section 8 in the paper [11] (see p. 171)
are a  special case of Theorem 5.1 (for the metric $\overline{\eta} \in {\rm WP}(g_M; a)$ with
$f(t) = ct$).

\bigskip \noindent
\section{The initial-value formulation for
the Einstein-Dirac equation }

\noindent
In this section we set up an invariant initial-value formulation
for the Einstein-Dirac equation     \[           {\rm Ric}_{\overline{\eta}} -
\frac{1}{2} S_{\overline{\eta}}  \overline{\eta}  =
T_{\overline{\eta}} ,  \qquad D_{\overline{\eta}} \psi =  \lambda_Q \psi,  \quad  \lambda_Q \in  {\mathbb R} ,        \]
where
\[         T_{\overline{\eta}} (X, Y) = \frac{\epsilon}{4} \Big( X \cdot \nabla^{\overline{\eta}}_Y
\psi +  Y \cdot \nabla^{\overline{\eta}}_X \psi, \ \psi \Big),   \qquad \epsilon = \pm 1 .
\]
The formulation will be applied in the next section to establish a local existence theorem
for a specific class of initial data sets(see Theorem 7.1).
Following the work [7] as a guideline, we can indeed express the evolution equations as well as the constraints in an invariant form.
For simplicity, we write  \[           \Delta : =  {\rm Ric}_{\overline{\eta}} -
\frac{1}{2} S_{\overline{\eta}}  \overline{\eta}  -
T_{\overline{\eta}} .       \]             The tensor field  $\Delta$ decomposes into three
parts      \[               \Delta = \Delta^E +   \Big\{  \Delta^M
\otimes F^{n+1}  +  F^{n+1}  \otimes  \Delta^M  \Big\} +
\Delta^H  \Big(  F^{n+1}  \otimes  F^{n+1}   \Big) ,      \]              where
\[
 \Delta^E   =   \sum_{i, j =1}^n   \Delta ( F_i, F_j ) F^i  \otimes F^j   ,    \quad
 \Delta^M   =   \sum_{i =1}^n   \Delta ( F_{n+1}, F_i )   F^i   ,   \quad
 \Delta^H   =    \Delta ( F_{n+1},  F_{n+1} )    .
\]

\bigskip \noindent
Restricting the equations,  $ \Delta^M = 0 $ and  $ \Delta^H  = 0$,
to a fixed slice, we obtain the {\it momentum constraint}    \[
T_{\overline{g}} ( F_{n+1}, V)    =   d( {\rm Tr}_{\overline{g}}
{\rm II} )(V) -  {\rm div}_{\overline{g}} ({\rm II} )(V)        \]       and
the {\it Hamiltonian  constraint}      \[  T_{\overline{g}}( F_{n+1},
F_{n+1} )   =   - \frac{1}{2}  S_{\overline{g}} +
\frac{1}{2} (  {\rm Tr}_{\overline{g}} {\rm II} )^2 -  \frac{1}{2}
{\rm Tr}_{\overline{g}} ( {\rm II}^2 ) ,         \]               where $
T_{\overline{g}}$ denotes the restriction of $
T_{\overline{\eta}}$ to the slices.
The information on the evolution should be contained in   \[
\Delta^E = 0  , \]          or any combination of it with the constraints.
The evolution equations should be chosen in such a way that, under the evolution, local preservation of the constraints
is guaranteed.
We consider the evolution equations of the form,     \[          \Delta (V, W)
=      \Delta ( F_{n+1}, F_{n+1} )   \cdot   \overline{\eta} (V,
W)  \quad  \mbox{and}   \quad   D_{\overline{\eta}} \psi = \lambda_Q \psi  ,   \quad  \lambda_Q \in {\mathbb R} . \]                                                      Note that, under
this evolution, the tensor $\Delta$ becomes         \[           \Delta  =
\Delta^M  \otimes F^{n+1}  +   F^{n+1}  \otimes  \Delta^M   +
\Delta^H  \cdot  \overline{\eta}  .       \]

\bigskip  \noindent
\begin{pro}    The
equation     \[          \Delta (V, W)  =      \Delta ( F_{n+1}, F_{n+1} )
\cdot   \overline{\eta} (V, W)         \]            is equivalent to
\begin{eqnarray*}
&      & {\rm Ric}_{\overline{\eta}} (V, W)          \\   &       &
\\
&   =   &  \frac{1}{n-1}   \Big\{  S_{\overline{g}} +
 {\rm Tr}_{\overline{g}} ( {\rm II}^2 ) - (  {\rm Tr}_{\overline{g}} {\rm II} )^2 -  {\rm Tr}_{\overline{\eta}} (T_{\overline{\eta}}) +
2  T_{\overline{\eta}}( F_{n+1}, F_{n+1} )   \Big\} \cdot  \overline{\eta} (V, W) +
T_{\overline{\eta}}(V, W)  .
\end{eqnarray*}
\end{pro}

\noindent
{\bf Proof.} Contracting the
 equation    \[       \Delta (V, W)  =      \Delta ( F_{n+1},
F_{n+1} )   \cdot   \overline{\eta} (V, W)  ,      \]              we see that
\[            \Delta ( F_{n+1}, F_{n+1} ) = - \frac{n-1}{2(n+1)} S_{\overline{\eta}} -
\frac{1}{n+1}  {\rm Tr}_{\overline{\eta}} (T_{\overline{\eta}})  ,            \]             from which it follows that                                     \[
{\rm Ric}_{\overline{\eta}} (V, W) =  \frac{1}{n+1}   \Big\{
S_{\overline{\eta}}
 -   {\rm Tr}_{\overline{\eta}} (T_{\overline{\eta}})    \Big\}  \cdot  \overline{\eta} (V, W)  + T_{\overline{\eta}}(V, W)  .         \]
Let us contract this equation. Then, with the help of the Gauss equation (iii) in Lemma 3.1, we
obtain         \[                   S_{\overline{\eta}} =  \frac{n+1}{n-1}  \Big\{
S_{\overline{g}} +
 {\rm Tr}_{\overline{g}} ( {\rm II}^2 ) - (  {\rm Tr}_{\overline{g}} {\rm II} )^2   \Big\}
- \frac{2}{n-1}  {\rm Tr}_{\overline{\eta}} (T_{\overline{\eta}}) +
\frac{2(n+1)}{n-1}   T_{\overline{\eta}}( F_{n+1}, F_{n+1} ) ,     \]                which gives the
asserted formula immediately. The converse is easy to verify.
$\rule{2mm}{2mm}$

\bigskip
Now we should verify that the constraints are indeed preserved under the
evolution
\[           \Delta (V, W)
=      \Delta ( F_{n+1}, F_{n+1} )   \cdot   \overline{\eta} (V,
W),  \qquad  D_{\overline{\eta}} \psi = \lambda_Q \psi .         \]                                                                                                     We note at this point that the divergence of the energy-momentum
tensor $T_{\overline{\eta}}$ vanishes, ${\rm div}_{\overline{\eta}} ( T_{\overline{\eta}} ) = 0$, since $T_{\overline{\eta}}$ is defined by eigenspinors of the Dirac operator $D_{\overline{\eta}}$ (see [12]).   Then, computing the
divergence, $ {\rm div}_{\overline{\eta}} (\Delta) = 0 $, and
expressing the covariant derivative $\nabla^{\overline{\eta}}$ in
terms of $\nabla^{\eta}$, we find that
\begin{eqnarray*}
0 & = & \sum_{j=1}^n  (d \Delta^H)( F_j ) F_j +  \sum_{j=1}^n  (
\nabla^{\eta}_{F_{n+1}} \Delta^M ) (F_j ) F_j
      \\
&      &      \\    &      & -  \Delta^H  \sum_{i,j=1}^n  \eta (
\Lambda_g (E_i, E_i),  \, E_j ) F_j    -  ( {\rm
Tr}_{\overline{g}} {\rm II} )   \sum_{j=1}^n  \Delta^M ( F_j ) F_j
\\ &      &       \\   &      &     -   \Delta^H  \sum_{i,j=1}^n  \eta
( E_i,  \, \Lambda_g (E_i, E_j) ) F_j        -
\sum_{j=1}^n  \Delta^M (  \rho^{-1} \nabla^g_{F_j}  \zeta +
\rho^{-1}  \Theta (F_j)  ) F_j
\end{eqnarray*}
and
\begin{eqnarray*}
0  & = &   \sum_{j=1}^n  ( \nabla^{\eta}_{F_j} \Delta^M ) (F_j)
F_{n+1} +  ( d \Delta^H ) ( F_{n+1} ) F_{n+1}   \\   &     &
\\    &     &  + \sum_{j=1}^n  \Delta^M ( ( \nabla^g_{F_j} L^{-1}
)(E_j) ) F_{n+1} + \sum_{j=1}^n  \Theta ( F_j, F_j )  \{  \rho
\Delta^H +  \Delta^M (\zeta)  \} F_{n+1}       \\   &     &     \\
&     & -  \sum_{j=1}^n  \Delta^M ( L^{-1} \{ \Lambda_{g} (E_j,
E_j)  \} ) F_{n+1}  + 2
\rho^{-1}  \sum_{j=1}^n d \rho (F_j )  \Delta^M (F_j) F_{n+1} ,
\end{eqnarray*}
where we have used the formula established in Proposition 2.1,     \[   \nabla^{\overline{g}}_{F_i} F_j  =
L^{-1} ( \nabla^g_{L^{-1} E_i} E_j ) +  L^{-1} \{ \Lambda_{g}
(E_i, E_j)  \} .           \]

\noindent Rewriting the above two equations with respect to $\eta$-orthonormal frame $( E_1, \ldots,
E_n, E_{n+1} )$, we arrive at
 a nonlinear hyperbolic system of first-order differential equations of the form
\[           \sum_{k=1}^n  A(k) \cdot  \nabla^{\eta}_{E_k} \Phi  +  B
\cdot  \nabla^{\eta}_{E_{n+1}} \Phi  + C \cdot \Phi    = 0 ,    \] where    \[
  \Phi   =    \Big\{  \sum_{j=1}^n   \Delta^M (E_j) E_j    \Big\} + (\Delta^H) E_{n+1} ,
\]
\[              A(k) =  \pmatrix{ - \rho^{-1}  \zeta^k (L^{-2})^i_j ,
&   \quad    (L^{-2})^i_k            \cr
   (L^{-2})^k_j,                                      &   \quad    - \rho^{-1} \zeta^k            \cr},     \qquad
B = \pmatrix{
 \rho^{-1}   (L^{-2})^i_j ,        &  \quad 0         \cr
0,         &    \quad   \rho^{-1}        \cr} ,          \]                                   and  $C$ is a (1,1)-tensor field given by
\begin{eqnarray*}
C^i_j & = & - ( {\rm Tr}_{\overline{g}} {\rm II} ) (L^{-2})^i_j - \rho^{-1}
\sum_{u=1}^n ( L^{-2} )^{iu} \eta ( \nabla^g_{E_u} \zeta + \Theta(E_u), \, E_j )  , \\
&     &     \\
C^i_{n+1} & = & - \sum_{u, v=1}^n  ( L^{-1} )^{iv}  \Big\{ \eta ( \Lambda_g ( E_u, E_u ), \, E_v)
+ \eta ( E_u, \, \Lambda_g ( E_u, E_v) ) \Big\}   ,   \\
&      &     \\
C^{n+1}_j & = &  \sum_{u, v=1}^n ( L^{-1} )^{uv} \eta ( ( \nabla^g_{E_u} L^{-1} )(E_v), \, E_j)
+ \sum_{u, v=1}^n \zeta^j ( L^{-2} )^{uv} \Theta (E_u, E_v)     \\
&      &     \\
&      & - \sum_{u=1}^n \eta ( L^{-1} \{ \Lambda_g ( E_u, E_u) \}, \, E_j ) +
\frac{2}{\rho} \sum_{u=1}^n ( L^{-2} )^u_j d \rho (E_u) ,     \\
&      &     \\
C^{n+1}_{n+1} & = & \rho  \sum_{u,v=1}^n ( L^{-2} )^{uv} \Theta (E_u, E_v) .
\end{eqnarray*}
We observe that the $(1,1)$-tensor fields $A(k)$
and $B$ are symmetric (with respect to reference metric $\eta$). Moreover, $B$ is
positive definite ($B \geq c I$ for some positive number $c > 0$), provided that every slice of $Q^{n+1}$ is compact. Thus, it is shown that, under our
evolution, the constraints are locally preserved.
Note that, when we consider the
warped product metrics as in Section 5 and Section 7, the local preservation of constraints holds without the assumption that every slice of $Q^{n+1}$ is compact.

\bigskip
Next, we
state a complete set of evolution equations for the Einstein-Dirac equation. Soon we will also define the corresponding initial data sets precisely.   Combining Corollary 3.1, Proposition 3.2, Proposition 4.1, 4.2 and Proposition 6.1 altogether, we easily obtain the evolution system of three differential equations, describing the
evolution of metrics $L^2 = \overline{g}$ , that of symmetric (0,2)-tensor fields ${\rm II}$  and that of spinor fields $\psi$ , respectively:
\begin{eqnarray*}
{\rm (E1)} &     & \eta ( ( \nabla^{\eta}_{E_{n+1}} L) (V), \ L(W)
) +  \eta (L(V), \ ( \nabla^{\eta}_{E_{n+1}} L) (W) )      \\   &
&     \\    & = & \eta ( ( \nabla^g_{\zeta} L)(V), \ L(W) ) + \eta (
L(V), \  ( \nabla^g_{\zeta} L)(W) )     \\   &     &     \\   &
& + \eta ( \nabla^g_V \zeta + \Theta(V) - \rho {\rm II} (V), \
L^2(W) ) +
 \eta ( L^2(V), \  \nabla^g_W \zeta + \Theta(W) - \rho {\rm II}(W) ) ,      \\
&       &       \\
&       &       \\
 {\rm (E2)}   &       & ( \nabla^{\eta}_{E_{n+1}} {\rm II}) (V, W)     \\ &
&      \\   &  =  &   \frac{\rho}{n-1}  \Big\{  S_{\overline{g}} +
{\rm Tr}_{\overline{g}} ( {\rm II}^2 ) -
  ( {\rm Tr}_{\overline{g}} {\rm II} )^2   \Big\}  \overline{g} (V, W)       \\
&       &            \\
&        & + \frac{\epsilon \rho}{4} \Big(   V \cdot
\nabla^{\overline{g}}_W \psi +
 W \cdot \nabla^{\overline{g}}_V \psi   ,   \  \psi   \Big)    \\
&       &      \\   &       & + \frac{\epsilon \rho}{8} \Big(   \{
V \cdot {\rm II} (W) + W \cdot {\rm II} (V) \} \cdot F_{n+1} \cdot
\psi , \  \psi  \Big)      \\   &       &      \\   &
& + \frac{\epsilon  \lambda_Q \rho}{2(n-1)}  ( \psi,  \   \psi
)  \overline{g} (V, W)  - \frac{\epsilon \rho}{n-1}
 \Big(  \sum_{i=1}^n F_i  \cdot  \nabla^{\overline{g}}_{F_i}  \psi ,  \  \psi    \Big) \overline{g} (V, W)       \\
&       &      \\   &       & +   \rho  \cdot  \Big\{  -  {\rm
Ric}_{\overline{g}}(V,W) - 2 \overline{g} (V,  {\rm II}^2(W) ) +
( {\rm Tr}_{\overline{g}} {\rm II} ) \cdot {\rm II}(V, W)   \Big\}
\\   &       &     \\   &       & + ( \nabla^{\overline{g}}_{\zeta}
{\rm II} ) (V,W) +  \overline{g} ( {\rm II}(V),
\nabla^{\overline{g}}_W \zeta ) +   \overline{g} ( {\rm II}(W),
\nabla^{\overline{g}}_V \zeta )         \\   &       &     \\   &
& + \overline{g} ( \Theta(V),  {\rm II}(W) )  +  \overline{g} (
\Theta(W),  {\rm II}(V) ) + (  \nabla^{\overline{g}}_W d \rho )
(V) ,           \\
&      &     \\
&      &      \\
{\rm (E3)}  &      &
 \nabla^{\eta}_{E_{n+1}} \psi      \\ &      &
\\   &  = &    \nabla^g_{\zeta} \psi  - \lambda_Q \rho
 E_{n+1} \cdot \psi     +  \frac{\rho}{2} (
{\rm Tr}_{\overline{g}} {\rm II} ) \psi          \\   &      &
\\  &       &  + \rho E_{n+1}  \cdot  \Big\{  \sum_{i=1}^n  E_i
\cdot   \nabla^g_{L^{-1} E_i} \psi +  \frac{1}{4}
\sum_{j,k,l=1}^n ( \Lambda_g )_{jkl} E_j  \cdot E_k  \cdot E_l
\cdot  \psi    \Big\}      \\   &       &       \\   &      &   +
\frac{1}{4} \sum_{i=1}^n E_i \cdot ( \nabla^{\eta}_{E_{n+1}} L )
(L^{-1} E_i) \cdot \psi + \frac{\rho}{4}   \sum_{i=1}^n  E_i
\cdot (L \circ {\rm II} \circ L^{-1})(E_i) \cdot \psi    \\  &
&     \\   &    &  - \frac{1}{4} \sum_{i=1}^n E_i \cdot (
\nabla^g_{\zeta} L ) (L^{-1} E_i) \cdot \psi      -
\frac{1}{4} \sum_{i=1}^n E_i \cdot L( \nabla^g_{L^{-1}E_i} \zeta )
\cdot \psi
   \\
&      &     \\   &      &   -   \frac{1}{4} \sum_{i=1}^n  E_i
\cdot (L \circ \Theta \circ L^{-1})(E_i) \cdot \psi    -
\frac{1}{2} \sum_{i=1}^n d \rho (L^{-1} E_i) E_i \cdot E_{n+1}
\cdot \psi.
\end{eqnarray*}

\bigskip  \noindent
In case of $n=2m-1$, the equation (E3) may be certainly replaced by the one in Corollary 4.1, and                                                                     the terms for ${\rm Tr}_{\overline{\eta}} (T_{\overline{\eta}})$
and $T_{\overline{\eta}} (F_{n+1}, F_{n+1})$ in (E2) may be replaced by the ones in Corollary 4.2.

\bigskip
Let us now define the initial data sets. We derive the constraint equations on initial hypersurfaces in a natural way,
by combining Proposition 4.2 (Corollary 4.2) with the relations,
\[       T_{\overline{\eta}} (V, F_{n+1}) = {\rm Ric}_{\overline{\eta}} (V, F_{n+1}) =
d({\rm Tr}_{\overline{g}} {\rm II} )(V) - {\rm
div}_{\overline{g}}({\rm II} )(V)        \]
and
\[         T_{\overline{\eta}} (F_{n+1}, F_{n+1}) = {\rm Ric}_{\overline{\eta}} (F_{n+1}, F_{n+1}) -
\frac{1}{2} S_{\overline{\eta}} = \frac{1}{2} \Big\{  - S_{\overline{g}}  + ( {\rm Tr}_{\overline{g}} {\rm II} )^2
- {\rm Tr}_{\overline{g}} (  {\rm II} \circ {\rm II} )   \Big\} .     \]

\bigskip  \noindent
{\bf Definition 6.1} (In case of $n=2m$)
   An {\it initial data set} $(
M^{2m}, \overline{g}, {\rm II}_M, \psi_M )$ for the
Einstein-Dirac equation on $Q^{2m+1}$ consists of a  slice
$M^{2m}$ with, defined on it, a metric $\overline{g}$,
a symmetric $(0,2)$-tensor field ${\rm II}_M$ and a spinor field
$\psi_M$  satisfying
the {\it momentum constraint}
\begin{eqnarray*}
&      &   d({\rm Tr}_{\overline{g}} {\rm II}_M )(V) - {\rm
div}_{\overline{g}}({\rm II}_M )(V)         \\   &      &     \\   &
=  &  \frac{\epsilon}{4}  \Big(  ( \sqrt{-1} )^{m+1}
\mu_{\overline{g}} \cdot  \Big\{ \nabla^{\overline{g}}_V  \psi_M
-  V \cdot D_{\overline{g}} \psi_M  \Big\}  , \      \psi_M
\Big),  \quad \epsilon = \pm 1,
\end{eqnarray*}
as well as the {\it Hamiltonian constraint}
\[
  - S_{\overline{g}}  + ( {\rm Tr}_{\overline{g}} {\rm II}_M )^2
- {\rm Tr}_{\overline{g}} (  {\rm II}_M \circ {\rm II}_M )
  =     - \epsilon \Big(
D_{\overline{g}}  \psi_M -  \lambda_Q
 \psi_M    ,        \   \psi_M      \Big),   \quad  \lambda_Q \in {\mathbb R}.
\]

\noindent
{\bf Definition 6.2} (In case of $n=2m-1$)
   An {\it initial data set} $(
M^{2m-1}, \overline{g}, {\rm II}_M, \psi^+_M, \varphi^+_M )$ for the
Einstein-Dirac equation on $Q^{2m}$ consists of a  slice
$M^{2m-1}$ with, defined on it, a metric $\overline{g}$,
a symmetric $(0,2)$-tensor field ${\rm II}_M$ and two spinor fields
$\psi^+_M, \varphi^+_M$  satisfying
the {\it momentum constraint}
\begin{eqnarray*}
&      &   d({\rm Tr}_{\overline{g}} {\rm II}_M )(V) - {\rm
div}_{\overline{g}}({\rm II}_M )(V)         \\ &      &     \\   &
=  &  \frac{\epsilon}{4}  \Big( \nabla^{\overline{g}}_V  \psi^+_M
-  V \cdot D_{\overline{g}} \psi^+_M    , \      \varphi^+_M \Big) -  \frac{\epsilon}{4}  \Big( \nabla^{\overline{g}}_V  \varphi^+_M
-  V \cdot D_{\overline{g}} \varphi^+_M    , \      \psi^+_M
\Big)
\end{eqnarray*}
as well as the {\it Hamiltonian constraint}
\begin{eqnarray*}
&   &  - S_{\overline{g}}  + ( {\rm Tr}_{\overline{g}} {\rm II}_M )^2
- {\rm Tr}_{\overline{g}} (  {\rm II}_M \circ {\rm II}_M )
\\   &     &     \\   &   =   &  - \epsilon \Big(
D_{\overline{g}}  \psi^+_M ,        \   \varphi^+_M     \Big) -  \epsilon \Big(
D_{\overline{g}}  \varphi^+_M ,        \   \psi^+_M     \Big) +
\epsilon \lambda_Q \Big\{ (\psi^+_M , \psi^+_M ) + ( \varphi^+_M, \varphi^+_M )  \Big\},  \quad  \lambda_Q \in {\mathbb R}.
\end{eqnarray*}

\bigskip  \noindent
\section{A local existence theorem }

\noindent
For a specific class of initial data sets, we will establish a local existence theorem for the Einstein-Dirac equation.
Let us begin with the case $n=2m$.
Let $\psi_M = \psi_M^+ + \psi^-_M$ be a spinor field on $(M^{2m}, g_M)$ with
 $\psi_M^{\pm} \in  \Gamma( \Sigma^{\pm} (M) )$.
Let $\Gamma_{\rm odd} ( \psi_M)$ denote the space of all spinor fields of the form
$\psi = h^+ \psi^+_M + h^- \psi^-_M$ on $Q^{2m+1} = M^{2m} \times {\mathbb R}$ defined by
\[      \psi(x, t) =   h^+(t) \psi^+_M(x) + h^-(t) \psi^-_M(x),   \qquad (x, t) \in M^{2m} \times {\mathbb R} ,        \]
where $h^{\pm} : {\mathbb R}  \longrightarrow {\mathbb C}$ are complex-valued functions.
The following lemma is an immediate consequence of Proposition 4.2, combined with Lemma 5.1.

\bigskip \noindent
\begin{lem}      Let $\psi_M = \psi^+_M + \psi^-_M$ be a real Killig spinor on $(M^{2m}, g_M)$ with
\[               \nabla^{g_M}_V \psi^{\pm}_M = - \frac{\lambda_M}{2m} V \cdot \psi^{\mp}_M,  \quad \lambda_M \in {\mathbb R} .           \]
Then $(\psi_M, \psi_M) = (\psi^+_M, \psi^+_M) + (\psi^-_M, \psi^-_M)$ is constant on $M^{2m}$, and the energy-momentum tensor, determined by $\overline{\eta} \in  {\rm WP}( g_M ; a)$   and $\psi \in  \Gamma_{\rm odd}( \psi_M)  $, is given by                                  \begin{eqnarray*}
{\rm Tr}_{\overline{\eta}} (T_{\overline{\eta}}) & = & \frac{\epsilon \lambda_Q}{2}
 \Big\{ ( h^+ \psi^+_M, h^+ \psi^+_M ) + ( h^- \psi^-_M, h^- \psi^-_M ) \Big\} ,
\\
&     &      \\
T_{\overline{\eta}} (V, W)  & = &  \frac{\epsilon \lambda_M}{4m} e^{\frac{f}{2}}
 \Big\{ ( h^+ \psi^-_M, h^- \psi^-_M ) + ( h^- \psi^+_M, h^+ \psi^+_M ) \Big\} \eta(V, W) ,
\\
&     &      \\
T_{\overline{\eta}} (V, F_{2m+1} ) & = & - \frac{ (2m+1) \epsilon \lambda_M}{8m}
 ( h^+ E_{2m+1} \cdot V \cdot \psi^-_M, \ h^- \psi^-_M )      \\
&      &      \\                                                                                &      &  - \frac{ (2m+1) \epsilon \lambda_M}{8m}
( h^- E_{2m+1} \cdot V \cdot \psi^+_M, \ h^+ \psi^+_M )  ,  \\
&     &      \\
T_{\overline{\eta}} ( F_{2m+1}, F_{2m+1} ) & = &
- \frac{\epsilon \lambda_M}{2} e^{- \frac{f}{2}}
 \Big\{ ( h^+ \psi^-_M, h^- \psi^-_M ) +  ( h^- \psi^+_M, h^+ \psi^+_M )  \Big\}       \\
&       &        \\                                                                            &       &   + \frac{\epsilon \lambda_Q}{2} \Big\{ ( h^+ \psi^+_M, h^+ \psi^+_M ) +  ( h^- \psi^-_M, h^- \psi^-_M )  \Big\} .
\end{eqnarray*}
\end{lem}

\bigskip \noindent
\begin{pro}      For $\overline{\eta} \in {\rm WP}( g_M ; a)$
and $\psi  \in \Gamma_{\rm odd} (\psi_M)$, the evolution equations (E1), (E2), (E3) for
the Einstein-Dirac equation are equivalent to
\begin{eqnarray*}
& (i) & \quad f_{tt}  =  \frac{a f_t f_t}{2}  - \frac{2}{m^2} ( \lambda_M )^2
e^{(a-1)f}   - \frac{\epsilon \lambda_Q}{2m-1} e^{af}
\langle h^+, h^+  \rangle  (\psi_M, \psi_M)    \\
&     &     \\
&     &   \quad \qquad   + \frac{2m+1}{4m(2m-1)} \epsilon \lambda_M e^{(a- \frac{1}{2})f}
\{ \langle h^+, h^-  \rangle  + \langle h^-, h^+  \rangle  \} (\psi_M, \psi_M) ,   \\
&     &     \\
&     &      \\
& (ii) & \quad h^+_t  =  - \frac{m}{2} f_t h^+ + ( \sqrt{-1} )^{2m+3} \lambda_Q e^{\frac{a}{2}f} h^+  -  ( \sqrt{-1} )^{2m+3}  \lambda_M e^{\frac{1}{2} (a-1)f} h^-  ,     \cr
&     &     \\
&     &  \quad  h^-_t   =   ( \sqrt{-1} )^{2m+3}  \lambda_M e^{\frac{1}{2} (a-1)f} h^+                                                                                         - \frac{m}{2} f_t h^- - ( \sqrt{-1} )^{2m+3} \lambda_Q e^{\frac{a}{2}f} h^- ,
\end{eqnarray*}
where we have used the notation $\langle \, , \, \rangle$ (for complex-valued
functions) to mean the standard Hermitian product.
\end{pro}

\noindent
{\bf Proof.} Proposition 4.1 implies that   $\psi = h^+ \psi^+_M + h^- \psi^-_M \in \Gamma_{\rm odd} (\psi_M)$ satisfies the Dirac equation,
\[        D_{\overline{\eta}} \psi  =  \lambda_Q  \psi,    \qquad   \overline{\eta} \in {\rm WP}( g_M ; a) ,       \]       on $(Q^{2m+1},  \overline{\eta})$ if and only if the second part (ii) of
the proposition is true. It remains to verify that the first part (i) of the proposition is
locally equivalent to the evolution equations (E1)-(E2).
Substituting Lemma 5.1, 5.2 and Lemma 7.1 in Proposition 6.1, we obtain
\begin{eqnarray*}
f_{tt} & = & \frac{a f_t f_t}{2}  - \frac{2}{m^2} ( \lambda_M )^2
e^{(a-1)f}  \\
&     &     \\                                                                                 &     & - \frac{\epsilon \lambda_Q}{2m-1} e^{af}
\Big\{ (h^+ \psi^+_M, h^+ \psi^+_M) + ( h^- \psi^-_M, h^- \psi^-_M )  \Big\}    \\
&     &     \\
&     & + \frac{2m+1}{2m(2m-1)} \epsilon \lambda_M e^{(a- \frac{1}{2})f}
\Big\{ ( h^+ \psi^-_M, h^- \psi^-_M) + ( h^- \psi^+_M, h^+ \psi^+_M )  \Big\} .
\end{eqnarray*}
Now we must note that the second part (ii) of the proposition implies
\[
\frac{d}{dt} \{ \langle h^+, h^+ \rangle - \langle h^-, h^- \rangle \}
 =  - m f_t  \{ \langle h^+, h^+ \rangle - \langle h^-, h^- \rangle \}   .
\]
Therefore, provided $h^+(0) = h^-(0)$ holds initially, the equality
$\langle h^+ , h^+ \rangle = \langle h^- , h^- \rangle$ is valid locally in $t$, and hence
\[            (h^+ \psi^+_M, h^+ \psi^+_M ) + ( h^- \psi^-_M, h^- \psi^-_M )
= \langle h^+, h^+ \rangle (\psi_M, \psi_M)         \]
is valid locally in $t$. Moreover, it is evident that
\[            ( h^+ \psi^-_M, h^- \psi^-_M ) + ( h^- \psi^+_M, h^+ \psi^+_M )  = \frac{1}{2} \{ \langle h^+ , h^- \rangle + \langle h^-, h^+  \rangle  \} ( \psi_M, \psi_M ).              \]
Thus we complete the proof of the proposition.    $\rule{2mm}{2mm}$

\bigskip \noindent
Let ${\rm Re}(h^{\pm})$ and ${\rm Im} (h^{\pm})$ denote the real and imaginary part of
the complex-valued functions $h^{\pm}$, respectively.
Then, we observe that, if we take       \[           \Psi = \Big( f   , \,  f_t , \, {\rm Re}(h^+),
\, {\rm Im}(h^+) ,  \, {\rm Re}(h^-) , \, {\rm Im}(h^-) \Big)        \]                        as a set of six unknowns, then
the system of evolution equations in Proposition 7.1 reduces to an autonomous eqution
\[          \frac{d}{dt} \Psi = H( \Psi )          \]
for some vector field $H$ defined on the six-dimensional Euclidean space ${\mathbb R}^6$.
This fact implies that, to each initial data, there corresponds a unique smooth local
solution to the evolution system in Proposition 7.1.

\bigskip \noindent
\begin{pro}
Let $(M^{2m}, g_M)$ be a Riemannian manifold admitting a real Killing spinor $\psi_M$. Then, for any real number $\lambda_Q \in {\mathbb R}$,
there exists an open interval $( - \omega, \omega) \subset {\mathbb R}$ and a warped product metric $\overline{\eta}$ on
$Q^{2m+1} = M^{2m} \times (- \omega, \omega)$ such that $(Q^{2m+1}, \overline{\eta})$ admits an Einstein spinor $\psi$ to
eigenvalue $\lambda_Q$.  In particular, if $\psi_M$ is a parallel spinor, then the Einstein spinor $\psi$ coincides with the WK-spinor in Proposition 5.2.
\end{pro}

\noindent
{\bf Proof.}  Let $\psi_M = \psi^+_M + \psi^-_M$ be a real Killing spinor, to Killing number $- \frac{\lambda_M}{2m} \in {\mathbb R}$, on the initial hypersurface $(M^{2m}, g_M)$. We identify $M^{2m}$ with the subspace
$M^{2m} \times \{ 0 \} \subset M^{2m} \times {\mathbb R}$.
Let  $\displaystyle \overline{\eta} = e^f \Big( \sum_{i=1}^{2m} E^i \otimes E^i \Big) +  e^{af}  dt \otimes
dt  \in {\rm WP}(g_M;a) $ and
$\psi = h^+ \psi^+_M + h^- \psi^-_M  \in \Gamma_{\rm odd} (\psi_M) $ satisfy the initial
conditions,
\[
 h^+(0) = h^-(0) = 1,   \qquad  f(0) = 0 ,    \] and
\[
 f_t(0) =  \pm  \sqrt{ \frac{4 (\lambda_M)^2}{m^2}
+ \frac{2 \epsilon (\lambda_Q - \lambda_M)}{m(2m-1)}  (\psi_M, \psi_M) } ,
\]
where we can always control $\epsilon = \pm 1$ and $(\psi_M, \psi_M) = \mbox{\it constant}$ so that
\[      \frac{4 (\lambda_M)^2}{m^2}
+ \frac{2 \epsilon (\lambda_Q - \lambda_M)}{m(2m-1)}  (\psi_M, \psi_M)               \]
is nonnegative. Let ${\rm II}_M = - \frac{1}{2} f_t(0) g_M$ be the symmetric (0,2)-tensor field required to prescribe initial data. Then, with the help of Lemma 7.1, one verifies that the initial data $( M^{2m}, g_M,$ $ {\rm II}_M, \psi_M)$ satisfies the constraint equations in Definition 6.1. Moreover, we know that the evolution system
in Proposition 7.1 is an autonomous equation and hence allows a local solution satisfying the initial data.  This proves the former part of the proposition.  Let us now suppose that the  spinor $\psi_M$ is a parallel spinor ($\lambda_M = 0$).  In this case,
we may assume that $\psi_M = \psi^+_M \in \Gamma (\Sigma^+ (M))$ and
$\psi = h^+ \psi^+_M$, and hence the evolution system in Proposition 7.1 simplifies to
\begin{eqnarray*}
f_{tt} & = & \frac{a f_t f_t}{2} - \frac{\epsilon \lambda_Q}{2m-1} e^{af}
\langle h^+, h^+ \rangle (\psi^+_M, \psi^+_M)   , \\
&     &     \\
h^+_t & = & - \frac{m}{2} f_t h^+ + ( \sqrt{-1} )^{2m+3} \lambda_Q e^{\frac{a}{2}f} h^+.
\end{eqnarray*}
On the other hand, since $\psi$ is (locally) an Einstein spinor, the
equation  \[            {\rm Ric}_{\overline{\eta}} (F_{2m+1}, F_{2m+1}) - \frac{1}{2} S_{\overline{\eta}}
= T_{\overline{\eta}} (F_{2m+1}, F_{2m+1})         \]
gives
\[                                                  \frac{m(2m-1)}{4} e^{-af} f_t f_t = \frac{\epsilon \lambda_Q}{2} \langle h^+, h^+  \rangle (\psi^+_M,  \psi^+_M).         \]
Thus, it follows that the function $f$ must satisfy $f_{tt} + \frac{m-a}{2} f_t f_t = 0$
whose solutions are exactly the ones given in Lemma 5.3, with $n=2m$.
$\rule{2mm}{2mm}$

\bigskip
Now, we proceed to the other case $n=2m-1$.
Let $\varphi^+_M$ be a spinor field on $M^{2m-1}$ (Note that we have identified $\Sigma(M)$ with
$\Sigma^+ (Q)$).  Let $\Gamma_{\rm even}( \varphi_M^+) $ denote the space of all spinor fields of the form $\varphi = h^+ \varphi^+_M +  k^+ E_{2m} \cdot \varphi^+_M$ on $Q^{2m} = M^{2m-1} \times {\mathbb R}$
defined by
\[       \varphi(x, t) =   h^+(t) \varphi^+_M(x) +  k^+(t) E_{2m} \cdot \varphi^+_M(x),   \qquad (x, t) \in M^{2m-1} \times {\mathbb R} ,         \]
where $h^+, \, k^+ : {\mathbb R}  \longrightarrow {\mathbb C}$ are complex-valued functions.
The following lemma is an immediate consequence of Corollary 4.2 combined with
Lemma 5.1.

\bigskip \noindent
\begin{lem}
Let $\varphi^+_M$ be a real Killig spinor on $(M^{2m-1}, g_M )$ with
\[          \nabla^{g_M}_V \varphi^+_M = - \frac{\lambda_M}{2m-1}  V \cdot E_{2m} \cdot \varphi^+_M.        \]
Then $(\varphi^+_M, \varphi^+_M )$ is constant on $M^{2m-1}$, and the energy-momentum tensor, determined by $ \overline{\eta} \in  {\rm WP}(g_M; a)$ and
$\varphi \in \Gamma_{\rm even}(\varphi^+_M)$, is given by
\begin{eqnarray*}
{\rm Tr}_{\overline{\eta}} (T_{\overline{\eta}})  & = &  \frac{\epsilon \lambda_Q}{2}
\Big\{ ( h^+ \varphi^+_M, h^+ \varphi^+_M ) + ( k^+ \varphi^+_M, k^+ \varphi^+_M )  \Big\} ,    \\
&     &    \\
T_{\overline{\eta}} (V, W) & = & \frac{\epsilon \lambda_M}{2m-1} e^{\frac{f}{2}}
( h^+ \varphi^+_M , k^+ \varphi^+_M ) \eta(V, W) ,        \\
&     &     \\
T_{\overline{\eta}} (V, F_{2m} ) & = & - \frac{m \epsilon \lambda_M}{2m-1} ( h^+ V \cdot E_{2m} \cdot \varphi^+_M , \
k^+ \varphi^+_M ) ,     \\
&     &     \\
 T_{\overline{\eta}} ( F_{2m} , F_{2m} ) & = &   - \epsilon \lambda_M e^{- \frac{f}{2}} ( h^+ \varphi^+_M , k^+ \varphi^+_M )
+  \frac{\epsilon \lambda_Q}{2}
\Big\{ ( h^+ \varphi^+_M, h^+ \varphi^+_M ) + ( k^+ \varphi^+_M, k^+ \varphi^+_M )  \Big\} .
\end{eqnarray*}
\end{lem}

\bigskip \noindent
\begin{pro}
For $\overline{\eta} \in {\rm WP}( g_M ; a)$
and $\varphi  \in \Gamma_{\rm even} (\varphi^+_M)$, the evolution equations (E1), (E2), (E3) for
the Einstein-Dirac equation are equivalent to
\begin{eqnarray*}
 & (i) &  \quad  f_{tt}
  =   \frac{a f_t f_t}{2} - \frac{8}{(2m-1)^2} ( \lambda_M )^2
e^{(a-1)f}    \\
&     &     \\
&     & \qquad  \qquad  - \frac{\epsilon \lambda_Q}{2(m-1)} e^{af}
\Big\{ (h^+ \varphi^+_M , h^+ \varphi^+_M ) + ( k^+ \varphi^+_M , k^+ \varphi^+_M ) \Big\}     \\
&     &     \\
&     &  \qquad  \qquad  + \frac{2m}{(m-1)(2m-1)} \epsilon \lambda_M e^{(a- \frac{1}{2})f}
(h^+ \varphi^+_M , k^+ \varphi^+_M )  ,     \\
&      &      \\
&      &       \\
 & (ii) &  \quad  h^+_t
 =  - \frac{2m-1}{4} f_t h^+ -  \lambda_M e^{\frac{1}{2}(a-1)f} h^+  +  \lambda_Q e^{\frac{a}{2}f} k^+  ,     \\
&     &     \\
&     &  \quad   k^+_t                                                                                           =  - \lambda_Q e^{\frac{a}{2} f} h^+                                                                                         - \frac{2m-1}{4} f_t k^+ + \lambda_M e^{\frac{1}{2}(a-1)f} k^+ ,
\end{eqnarray*}
where we may choose $h^+, k^+$ to be real-valued functions.
\end{pro}

\noindent
{\bf Proof.}  Corollary 4.1 implies that $\varphi = h^+ \varphi^+_M + k^+ E_{2m} \cdot \varphi^+_M \in \Gamma_{\rm even} (\varphi^+_M)$ satisfies the Dirac equation,
\[        D_{\overline{\eta}} \varphi =   \lambda_Q \varphi ,    \qquad  \overline{\eta} \in {\rm WP}(g_M; a) ,       \]
on $(Q^{2m}, \overline{\eta} )$ if and only if the second part (ii) of the proposition is true, where $h^+, k^+$ may be chosen to be real-valued functions.    Substituting Lemma 5.1, 5.2 and 7.2 in Proposition 6.1, we obtain the first part (i)
of the proposition.
$\rule{2mm}{2mm}$

\bigskip \noindent
With the help of Proposition 7.2 and 7.3, we prove the main theorem of the paper.

\bigskip \noindent
\begin{thm}
Let $(M^n, g_M)$ be a Riemannian manifold admitting a real Killing spinor $\varphi_M$. Then, for any real number $\lambda_Q \in {\mathbb R}$,
there exists an open interval $( - \omega, \omega ) \subset {\mathbb R}$ and a warped product metric $\overline{\eta}$ on
$Q^{n+1} = M^{n} \times ( - \omega, \omega)$ such that $(Q^{n+1}, \overline{\eta})$ admits an Einstein spinor $\varphi$ to
eigenvalue $\lambda_Q$.  In particular, if $\varphi_M$ is a parallel spinor, then the
Einstein spinor $\varphi$ coincides with the WK-spinor in Theorem 5.1.
\end{thm}

\noindent
{\bf Proof.}  Because of Proposition 7.2, it suffices to prove the theorem for the case $n=2m-1$. Let $\varphi^+_M $ be a real Killig spinor, to Killing number $- \frac{\lambda_M}{2m-1} \in {\mathbb R}$, on the initial hypersurface $(M^{2m-1}, g_M)$.
Let  $\displaystyle \overline{\eta} = e^f \Big( \sum_{i=1}^{2m-1} E^i \otimes E^i \Big) +  e^{af} dt \otimes
dt   \in {\rm WP}(g_M;a) $ and
$\varphi = h^+ \varphi^+_M + k^+ E_{2m} \cdot \varphi^+_M  \in \Gamma_{\rm even} (\varphi^+_M) $ satisfy the initial
conditions,
\[
h^+(0) = k^+(0) = 1,   \qquad  f(0) = 0 ,        \]
and
\[             f_t(0) =  \pm  \sqrt{ \frac{16 (\lambda_M)^2}{(2m-1)^2}
+ \frac{4 \epsilon ( \lambda_Q - \lambda_M) }{(m-1)(2m-1)} (\varphi^+_M, \varphi^+_M) },
\]
where we can always control $\epsilon = \pm 1$ and $(\varphi_M, \varphi_M) = \mbox{\it constant}$ so that
\[           \frac{16 (\lambda_M)^2}{(2m-1)^2}
+ \frac{4 \epsilon ( \lambda_Q - \lambda_M) }{(m-1)(2m-1)} (\varphi^+_M, \varphi^+_M)
  \]
is nonnegative. Let ${\rm II}_M = - \frac{1}{2} f_t(0) g_M$ be the symmetric (0,2)-tensor field required to prescribe initial data. Then, with the help of Lemma 7.2, one verifies that the initial data set $( M^{2m}, g_M,$ $ {\rm II}_M, \psi^+_M = \varphi^+_M, \varphi^+_M)$ satisfies the constraint equations in Definition 6.2. Moreover, as in the case of $n = 2m$, we find that there exists a unique local solution to the evolution system in Proposition 7.3 satisfying the initial data.  One proves the latter part
of the theorem in a similar way as in the proof for Proposition 7.2.   \
$\rule{2mm}{2mm}$

\bigskip \noindent
The Einstein spinors of Theorem 7.1 do not generally extend to
$M^n \times {\mathbb R} $, since the evolution system in Proposition 7.1 (resp. Proposition 7.3), in general,
do not allow global solutions. However, via reparametrization $(- \omega, \omega)
\longrightarrow {\mathbb R}$, we conclude that, indeed, there exist global solutions
to the Einstein-Dirac equation on $M^n \times {\mathbb R}$.

\bigskip \noindent
\begin{cor}
Let $(M^n, g_M)$ be a Riemannian manifold admitting a real Killing spinor. Then, for any real number $\lambda_Q \in {\mathbb R}$,
there exists a warped product  metric $\overline{\eta}^{\ast}$ on
$Q^{n+1} = M^{n} \times {\mathbb R}$ such that $(Q^{n+1}, \overline{\eta}^{\ast} )$ admits an Einstein spinor to eigenvalue $\lambda_Q$.
\end{cor}

\noindent
{\bf Proof.} We consider the case $n=2m$. The same argument is valid for the other case $n=2m-1$. By Theorem 7.1, there exists a solution $( \overline{\eta}, \psi )$ to the Einstein-Dirac equation on $M^{2m} \times ( - \omega, \omega)$ for some positive number $\omega$,
with $\displaystyle \overline{\eta} = e^f \Big( \sum_{i=1}^{2m} E^i \otimes E^i \Big) +  e^{af} dt \otimes
dt   \in {\rm WP}(g_M;a) $ and
 $\psi = h^+ \psi^+_M + h^- \psi^-_M  \in \Gamma_{\rm odd} (\psi_M) $.
Let $\gamma : {\mathbb R}  \longrightarrow
( - \omega, \omega)$ be a diffeomorphism, e.g., defined by
\[    \gamma(s) = \frac{2 \omega}{\pi} {\rm arctan}(s), \qquad s \in {\mathbb R} .    \]
Now we pull back the metric $\overline{\eta}$ as well as the Einstein spinor $\psi$
to $M^{2m} \times {\mathbb R}$
via the diffeomorphism
\[       I \times \gamma : M^{2m} \times {\mathbb R} \longrightarrow M^{2m} \times ( - \omega, \omega) ,
\qquad (x, s ) \longmapsto (x, \gamma(s) ) .        \]

\noindent
In fact, using the diffeomorphism $\gamma$ and the relations,
\[      \frac{dt}{ds} = \frac{2 \omega}{\pi (s^2 + 1)} ,  \qquad
\frac{ds}{dt} = \frac{\pi}{2 \omega} (s^2 + 1) ,      \]
we can explicitly express the pullbacked objects as
\[
\overline{\eta}^{\ast}  : =  (I \times \gamma)^{\ast}(\overline{\eta})
= e^{f^{\ast}} \Big( \sum_{i=1}^{2m} E^i \otimes E^i \Big) +  e^{a {f^{\ast}}}
\frac{4 \ \omega^2}{\pi^2 (s^2+1)^2}  ds \otimes ds ,      \]
where $f^{\ast}(s) = (f \circ \gamma)(s) $,
and
\[             \psi^{\ast} : = (I \times \gamma)^{\ast} (\psi) =
(h^+ \circ \gamma) \psi^+_M +
(h^- \circ \gamma)  \psi^-_M .        \]
Obviously, $(\overline{\eta}^{\ast} , \psi^{\ast} )$ is a global solution to
the Einstein-Dirac equation on $M^{2m} \times {\mathbb R}$.    
$\rule{2mm}{2mm}$

\bigskip \bigskip  \noindent
{\bf Acknowledgements:} The author thanks Thomas Friedrich for reading a preliminary version
of the paper. The Section 5 of the paper was inspired by his comments. This research was
supported by the BK 21 project of Seoul National University, the BK 21 project of Inha
University and in part by the SFB 288 of the Deutsche Forschungsgemeinschaft.

\bigskip \noindent

\begin{Literature}{xx}
\bibitem{1}
D. Bao, J. Isenberg and P.B. Yasskin, {\it The dynamics of the
Einstein-Dirac system ; A principal bundle formulation of the
theory and its canonical analysis,} Annals of Physics 164 (1985) 103-171.
\bibitem{2}
B\"{a}r, {\it Real Killing spinors and holonomy,} Commum. Math. Phys. 154 (1993) 509-521.
\bibitem{3}
H. Baum, Th. Friedrich, R. Grunewald, I. Kath, {\it Twistors and Killing spinors on Riemannian
manifolds,} Teubner, Leipzig/Stuttgart (1991).
\bibitem{4}
F. Belgun, {\it The Einstein-Dirac equation on Sasakian 3-manifolds,} SFB 288 No. 477, Berlin (2000).
\bibitem{5}
J.P. Bourguignon, P. Gauduchon, Spineurs, {\it Op\'{e}rateurs de Dirac et Variations de
M\'{e} triques,} Commum. Math. Phys. 144 (1992) 581-599.
\bibitem{6}
Y. Choquet-Bruhat, {\it The Cauchy problem.} In: A. Held (Ed.) General Relativity and Gravitation,
Vol. I. Plenum, New York (1980) 99-172.
\bibitem{7}
H. Friedrich and A. Rendall, {\it The Cauchy problem for the Einstein
equations,} gr-qc/0002074.
\bibitem{8}
Th. Friedrich, {\it Der erste Eigenwert des Dirac-Operators einer kompakten Riemannschen
Mannigfaltigkeit nichtnegativer Skalarkr\"{u}mmung,} Math. Nachr. 97 (1980) 117-146.
\bibitem {9} Th. Friedrich, {\it Dirac operators in Riemannian geometry,} Graduate Studies in Math. 25, AMS (2000).
\bibitem{10} Th. Friedrich, {\it Solutions of the Einstein-Dirac equation on Riemannian 3-manifolds
with constant scalar curvature,} Journ. Geom. Phys. 36 (2000) 199-210.
\bibitem{11}
Th. Friedrich, E.C. Kim, {\it The Einstein-Dirac equation on Riemannian spin manifolds,}
Journ. Geom. Phys. 33 (2000) 128-172.
\bibitem{12}
E.C. Kim, {\it Die Einstein-Dirac-Gleichung \"{u}ber Riemannschen Spin-Mannigfaltig-keiten,}
Dissertation, Humboldt-Univ., Berlin (1999).
\bibitem{13}
P. Ramacher, {\it Geometric and analytic properties of families of hypersurfaces in Eguchi-Hanson
space,} SFB 288 No. 500 (2000).
\bibitem{14}
P. Tondeur, {\it Foliations on Riemannian manifolds,} Universitext, Springer-Verlag (1988).
\bibitem{15}
M. Wang, {\it Parallel spinors and parallel forms,} Ann. Glob. Anal. Geom. 7 (1989) 59-68.
\end{Literature}

\end{document}